\documentclass[12pt]{article}
\usepackage[top=1in,bottom=1.5in,left=1.1in,right=1.1in]{geometry}
\usepackage{amsmath,amssymb,amsfonts, euscript}
\usepackage{graphicx}
\DeclareGraphicsExtensions{.jpg,.png,.pdf,.ps} 
\usepackage{float}
\newtheorem{theorem}{Theorem}[section]
\newtheorem{proposition}[theorem]{Proposition}
\newtheorem{lemma}[theorem]{Lemma}
\newtheorem{corollary}[theorem]{Corollary}
\newtheorem{example}{Example}[section]
\newtheorem{definition}{Definition}[section]
\newtheorem{remark}{Remark}[section]
\newcommand{\repeatthanks}{\textsuperscript{\thefootnote}}
\title{Continuous-time GARCH process driven by semi-L\'evy process}
\author{M. Mohammadi \thanks{Faculty of Mathematics and Computer Science, Amirkabir University of Technology, 424 Hafez
Avenue, Tehran 15914, Iran. E-mail: m.mohammadiche@aut.ac.ir(M. Mohammadi) and rezakhah@aut.ac.ir(S. Rezakhah).}
\qquad S. Rezakhah\repeatthanks
\qquad N. Modarresi\thanks{Department of Mathematics and computer science, Allameh Tabataba’i University, Tehran, Iran. E-mail: n.modarresi@atu.ac.ir(N. Modarresi).}}
\begin{document}
\maketitle
\begin{abstract}
In this paper we study the simple semi-L\'evy driven continuous-time generalized autoregressive conditionally
heteroscedastic (SS-COGARCH) process. The statistical properties of this process are characterized. This process
has the potential to approximate any semi-L\'evy driven COGARCH processes. We show that the state representation
of such SS-COGARCH process can be described by a random recurrence equation with periodic random coefficients. 
The almost sure absolute convergence of the state process is proved. The periodically stationary solution of the
state process is shown which cause the volatility to be periodically stationary under some suitable conditions.
Also it is shown that the increments with constant length of such SS-COGARCH process is itself a periodically 
correlated (PC) process. Finally, we apply some test to investigate the PC behavior of the increments (with
constant length) of the simulated samples of proposed SS-COGARCH process.

\end{abstract}

\section{Introduction}
Many financial data and indices have  heteroscedastic structure. Examples of this kind are stocks returns, 
network traffic and natural data, see \cite{b5, k2, j}. Popular model for these data are autoregressive 
conditionally heteroscedastic (ARCH) model proposed by Engle \cite{e1} and generalized ARCH (GARCH),
Bollerslev \cite{b4}. The GARCH type processes have become the most popular tools to model heteroscedasticity in
discrete time.\\
In practice, for various reasons such as high-frequency data, many time series are irregularly spaced and this
has created a demand for continuous-time models, \cite{b9}. For the first time, Kluppelberg et al. \cite{k1}
introduced a continuous-time version of the GARCH(1,1) (COGARCH(1,1)) process, which preserves the
essential features of the discrete-time GARCH(1,1) processes. They replaced the noise of the discrete-time 
GARCH(1,1) process with the increments of some L\'evy process. The volatility of this process satisfies a
stochastic differential equation. They proved the stationarity property and also second order properties under
some regularity conditions on the corresponding L\'evy process. Brockwell et al. \cite{b9} generalized the 
L\'evy driven COGARCH(1,1) process to the L\'evy driven COGARCH$(p,q)$ process for $q\geq p\geq1$ when its 
volatility is a continuous-time ARMA (CARMA) process \cite{b8}. They showed that the state representation of the
volatility can be expressed as a stochastic recurrence equation with random coefficients.\\ 
Periodic behavior is common in many real-world time series such as power market prices, car accident claims for 
an insurance company and sales with seasonal interest. The term periodically correlated (PC) was introduced by 
Gladyshev \cite{g1}, but the same property was introduced by Bennett \cite{b2} who called them cyclostationary 
(\cite{h1}). Properties of PC processes are studies by Hurd and Miamee \cite{h1}. Bibi and Lescheb \cite{b3} 
studied the class of bilinear processes with periodic time-varying coefficients of periodic ARMA and periodic 
GARCH models.\\ 
L\'evy processes introduced by L\'evy have stationary and independent increments and right continuous paths with
left limits \cite{s1}. Such processes have potential to be applied to financial data following stochastic
volatility structure. A generalization of L\'evy process is semi-L\'evy process, that has periodically 
stationary increments, studied by Maejima and sato \cite{m1}. We considered this process as the underlying 
process in CARMA \cite{b8} and COGARCH \cite{b9, k1} processes that can be applied when there is evident that 
the underlying process has PC increments. The observations of such processes have significant dependency to the 
ones of previous periods. So semi-L\'evy process are more prominent than L\'evy processes in such cases.\\
In this paper we introduce a COGARCH process driven by some simple semi-L\'evy process, which we call SS-COGARCH
process. The simple semi-Levy process is defined as a compound Poisson process with periodic time-varying 
intensity with period $\tau$. This process enables us to provide the statistical properties of the SS-COGARCH 
process. Moreover, we find a random recurrence equation with periodic random coefficients for the state 
representation of such process. By some regularity condition we show the absolute convergence of the state 
equation. We also show that the volatility of the SS-COGARCH process is strictly periodically stationary. The 
increments of the SS-COGARCH process with constant length $h=\tau/\varrho$ where $\varrho$ is some integer is a 
discrete-time PC process with period $\varrho$. Such SS-COGARCH process has the potential to provide an
approximation for every semi-L\'evy driven COGARCH process. Finally, we investigate the theoretical results 
concerning PC structure of the increment process by simulation. We show that the increments of the SS-COGARCH 
process with length $h$ is PC with some period $\varrho$ and the support of the squared coherence statistics 
consists of lines parallel to the main diagonal and having spacing of $2\pi/\varrho$. \\
This paper is organized as follows. In section 2 we introduce the simple semi-Levy driven COGARCH processes. For
this, we present the simple semi-Levy process and obtain the characteristic function of it. Section 3 is devoted
to some sufficient conditions which make the volatility process strictly periodically stationary. We obtain the 
mean, covariance function of the state process and volatility process in section 4. We also investigate second 
order properties of the squared increments of the COGARCH process in this section. In section 5 we illustrate 
the results with simulations. All proofs are contained in Section 6.
 
\renewcommand{\theequation}{\arabic{section}.\arabic{equation}} 
\section{Simple semi-L\'evy driven COGARCH processes}
\setcounter{equation}{0}
In this section we study the preliminaries such as the additive processes and their characteristic functions and
semi-L\'evy process in subsection 2.1. We also describe the structure of simple semi-L\'evy process and 
characteristics it in subsection 2.2. Then we introduce the simple semi-L\'evy driven COGARCH (SS-COGARCH)
process in subsection 2.3.
\subsection{Preliminaries}
Let $(\Omega, \mathcal{F}, (\mathcal{F}_{t})_{t\geq0}, \mathbb{P})$ be a filtered probability space, where 
$\mathcal{F}_{t}$ is the smallest right-continuous filtration such that $\mathcal{F}_{0}$ contains all the
$\mathbb{P}$-null sets of $\mathcal{F}$. A process $(X_{t})_{t\geq0}$ defined on the probability space $(\Omega,
\mathcal{F}, (\mathcal{F}_{t})_{t\geq0}, \mathbb{P})$ is called an additive process if $X_{0}=0$ a.s., it is
stochastically continuous, it has independent increments and its sample paths are right-continuous and have left
limits in $t > 0$. Further, if $(X_{t})_{t\geq0}$ has stationary increments, it is a L\'evy process \cite{c2,
s1}. The characteristic function of the additive process $(X_{t})_{t\geq0}$ has a following L\'evy-Khinchin
representation \cite[Theorems 9.1-9.8]{s1}.
\begin{theorem}
Let $(X_{t})_{t\geq0}$ be an additive process on $\mathbb{R}^{d}$. Then $(X_{t})_{t\geq0}$ has infinitely 
divisible distribution for $t\geq0.$ The law of $(X_{t})_{t\geq0}$ is uniquely determined by its spot 
characteristic triplet $(\Gamma_{t},\Pi_{t},\psi_{t})_{t\geq0}$
\begin{align*}
E[e^{i<w,X_{t}>}]=e^{\varphi_{t}(w)},\qquad w\in\mathbb{R}^{d},
\end{align*}
\begin{align*}
\varphi_{t}(w)=i<w,\Gamma_{t}>-\frac{1}{2}<w,\Pi_{t}w>+\int_{\mathbb{R}^{d}}(e^{i<w,x>}-1-i<w,x>I_{\lbrace||x||
\leq1\rbrace})\psi_{t}(dx).
\end{align*}
\end{theorem}
where $<\cdot,\cdot>$ is inner product and $||\cdot||$ is Euclidean vector norm. The spot L\'evy measure $\psi_{t}$
satisfies the integrability condition 
$\int_{\mathbb{R}^{d}}min\lbrace 1, ||x||^{2}\rbrace \psi_{t}(dx)<\infty$ for $t\geq0.$
\begin{remark} By \cite[p.458-459]{c2}, the spot characteristic triplet $(\Gamma_{t}, \Pi_{t}, 
\psi_{t})_{t\in[0,T]}$ can be defined by
\begin{align*}
\Gamma_{t}&=\int_{0}^{t}\gamma_{s}ds\\
\Pi_{t}&=\int_{0}^{t}\sigma_{s}^{2}ds\\
\psi_{t}(B)&=\int_{0}^{t}\upsilon_{s}(B)ds,\qquad \forall B\in\mathcal{R}^{d},
\end{align*}
where $\mathcal{R}^{d}$ is $\sigma-$field on the $\mathbb{R}^{d}$. The triplet $(\gamma_{t}, \sigma_{t}^{2}, 
\upsilon_{t}) _{t\in[0,T]}$ is called the local characteristic triplet of $(X_{t})_{0\leq t\leq T}$ which satisfy the
following conditions:
\begin{itemize}
\item
$\gamma_{t}: [0, T]\rightarrow\mathbb{R}^{d}$ is a deterministic function with finite variation.
\item
$\sigma_{t}: [0,T]\rightarrow M_{d\times d}(\mathbb{R})$ is a symmetric, continuous and matrix valued function 
which verifies
$\int_{0}^{T}\sigma^{2}_{t}dt<\infty$.
\item
$(\upsilon_{t})_{t\in[0,T]}$ is a family of L\'evy measures which verifies
\begin{align*}
\int_{0}^{T}\int_{\mathbb{R}^{d}}min\lbrace1, ||x||^{2}\rbrace\upsilon_{t}(dx)dt<\infty.
\end{align*}
\end{itemize}
\end{remark}
As an extension of L\'evy process, we present the definition of semi-L\'evy processes \cite{m1}.
\begin{definition}
A subclass of additive processes $(X_{t})_{t\geq0}$ is called semi-L\'evy process with period $\tau>0,$ if for 
any $0\leq s\leq t$,
\begin{align*}
X_{t}-X_{s}\overset{d}{=}X_{t+\tau}-X_{s+\tau}
\end{align*}
where $\overset{d}{=}$ denotes the equality in distributions.
\end{definition}
\subsection{Structure of simple semi-L\'evy process}
For describing the structure of the simple semi-L\'evy process, we define the general structure of the 
intensities function of the Poisson process with periodically stationary increments. We also characterize this 
pure jump process by representation the characteristic function and introduce the corresponding semi-L\'evy measure.
\begin{definition}: \textbf{Poisson process with periodically stationary increment}\\
A process $\big(N(t)\big)_{t\geq0}$ is a Poisson process with periodically stationary increment where 
$E\big(N(t)\big)$ $=\Lambda(t),$
\begin{align}\label{equ2.1}
\Lambda(t)=\int_{0}^{t}\lambda(u)du
\end{align}
and the intensity $\lambda(\cdot)$ is a periodic non-negative function with some period $\tau>0,$ so 
$\lambda(t)=$ $\lambda(t+k\tau)$ for $t\geq0, k\in\mathbb{N}.$
\end{definition}
\begin{definition}: \textbf{Simple compound Poisson process}\\
Let $0=t_{0}<t_{1}<\cdots$ be a partition of positive real line. Also assume that $A_{j}=[t_{j-1},t_{j}),$ 
$j\in\mathbb{N},$ and $|A_{j}|=|A_{j+l}|$ for some integer $l\in\mathbb{N}$ and $\tau=\sum_{j=1}^{l}|A_{j}|.$ Let
$\big(N(t)\big)_{t\geq0}$ be a Poisson process which has periodically stationary increments with period $\tau>0$ and
intensity function $\Lambda(t)$ defined by (\ref{equ2.1}). Then the simple compound Poisson process $(S_{t})_{t\geq0}$
is defined as  
\begin{align}\label{equ2.2}
S_{t}=D_{t}+\sum_{n=1}^{N(t)}Z_{n}
\end{align}
where $Z_{n}=\sum_{j=1}^{l}Z_{n}^{j}I_{\lbrace \Upsilon_{n}\in\mathfrak{D}_{j}\rbrace},$ $\Upsilon_{n}$ is the 
arrival time of $n^{th}$ jump $Z_{n}$, $\mathfrak{D}_{j}=\bigcup_{k=0}^{\infty}A_{j+kl}$ and $Z_{n}^{j}$ are 
independent and have distribution $F_{j}, j=1, \cdots,l,$ such that $\int_{\mathbb{R}}z^{2}F_{j}(dz)<\infty$ for $j=1,
\cdots,l.$ Also $D_{t},$ $t>0,$ is a deterministic drift function with period $\tau$, say $D_{t}=D_{t+\tau},$ and 
$D_{0}=0.$ One can easily verify that $(S_{t})_{t\geq0}$ has independent increment.
\end{definition}
Now we find characteristic function of the simple compound Poisson process $(S_{t})_{t\geq0}$ by the following 
Lemma.
\begin{lemma}
Let $(N(t))_{ t\geq0}$ be a Poisson process with periodically stationary increment and mean $\Lambda(t)$, 
defined by (\ref{equ2.1}). Then the
process $(S_{t})_{t\geq0}$ defined by (\ref{equ2.2}) has the following characteristic function for $t\geq0$
\begin{align*}
E[e^{iwS_{t}}]=e^{\varphi_{t}(w)},
\end{align*}
\begin{align*}
\varphi_{t}(w)=iw\Gamma_{t}+\int_{\mathbb{R}}(e^{iw\mathit{z}}-1-iw\mathit{z} I_{\lbrace|\mathit{z}|
\leq1\rbrace})\psi_{t}(d\mathit{z}),
\end{align*}
where
\begin{align*}
\Gamma_{t}&=D_{t}+\sum_{k=0}^{m-1}\sum_{r=1}^{l}\int_{|z|
\leq1}z\big(\Lambda(t_{kl+r})-\Lambda(t_{kl+r-1})\big)F_{r}(dz)\\
&\quad+\sum_{r=1}^{j-1}\int_{|z|\leq1}z\big(\Lambda(t_{ml+r})-\Lambda(t_{ml+r-1})\big)F_{r}(dz)\\
&\quad+\int_{|z|\leq1}z\big(\Lambda(t)-\Lambda(t_{ml+j-1})\big)F_{j}(dz),
\end{align*}
and
\begin{align*}
\psi_{t}(dz)&=\sum_{k=0}^{m-1}\sum_{r=1}^{l}\big(\Lambda(t_{kl+r})-\Lambda(t_{kl+r-1})\big)F_{r}(dz)\\
&\quad+\sum_{r=1}^{j-1}\big(\Lambda(t_{ml+r})-\Lambda(t_{ml+r-1})\big)F_{r}(dz)\\
&\quad+\big(\Lambda(t)-\Lambda(t_{ml+j-1})\big)F_{j}(dz),
\end{align*}
where $m=[\frac{t}{\tau}]$ and $(t-m\tau)\in A_{j}$ for some $j=1,\cdots,l.$
\end{lemma}
Proof: see Appendix, P1.
\begin{remark}
By Remark 2.1, Lemma 2.2 and (2.1), the spot characteristic triplet of process $(S_{t})_{0\leq t\leq T},$ 
$(\Gamma_{t},0,\psi_{t})_{t\in[0,T]},
$ have the local characteristic triplet $(\gamma_{s}, 0, \upsilon_{s})_{s\in[0,T]}$ which has the following form
\begin{align*}
\gamma_{s}&=dD_{s}+\sum_{k=0}^{m-1}\sum_{r=1}^{l}\int_{|z|\leq1}z\lambda(s)I_{[t_{kl+r-1},t_{kl+r})}
(s)F_{r}(dz)\\
&\quad+\sum_{r=1}^{j-1}\int_{|z|\leq1}z\lambda(s)I_{[t_{ml+r-1},t_{ml+r})}(s)F_{r}(dz)\\
&\quad+\int_{|z|\leq1}z\lambda(s)I_{[t_{ml+j-1},t]}(s)F_{j}(dz),
\end{align*}
and
\begin{align}\nonumber
\upsilon_{s}(dz)&=\sum_{k=0}^{m-1}\sum_{r=1}^{l}\lambda(s)I_{[t_{kl+r-1},t_{kl+r})}(s)F_{r}(dz)\\ \nonumber
&\quad+\sum_{r=1}^{j-1}\lambda(s)I_{[t_{ml+r-1},t_{ml+r})}(s)F_{r}(dz)\\ \label{equ2.3}
&\quad+\lambda(s)I_{[t_{ml+j-1},t]}(s)F_{j}(dz).
\end{align}
\end{remark}
It follows from definition 2.3 and Remark 2.2 that the family $(\upsilon_{s})_{s\in[0,T]}$ of semi-L\'evy measures verify
\begin{align}\label{equ2.4}
\int_{0}^{T}\int_{\mathbb{R}}|z|^{2}\upsilon_{s}(dz)ds<\infty.
\end{align}
This implies that $(S_{t})_{t\geq0}$ is semi-martingale, so it has L\'evy-Ito decomposition and has
quadratic variation process \cite[p.459-460]{c2}. 
\begin{corollary}
By lemma 2.2, the stochastic process $(S_{t})_{t\geq0}$ defined by (\ref{equ2.2}) is a semi-L\'evy process with 
period $\tau.$
\end{corollary}
Proof: see Appendix, P2.

\subsection{Structure of simple semi-L\'evy driven COGARCH process}
Let $(S_{t})_{t\geq0}$ be a simple semi-L\'evy process  with period $\tau$ defined by (\ref{equ2.2}). Process
$(G_{t})_{t\geq0}$ with parameters $\alpha_{0}>0, \alpha_{1}, \cdots, \alpha_{p}\in\mathbb{R}, \beta_{1}, \cdots, 
\beta_{q}\in\mathbb{R}, \alpha_{p}\neq0, \beta_{q}\neq0,$ and $\alpha_{p+1}=\cdots\alpha_{q}=0$ is a simple semi-
L\'evy driven COGARCH(p,q) process (SS-COGARCH(p,q)), $q\geq p\geq1,$ defined by $dG_{t}=\sqrt{V_{t}}dS_{t}$ or
equivalently 
\begin{equation}\label{equ2.5}
G_{t}=\int_{0}^{t}\sqrt{V_{u}}dS_{u},\qquad t>0,\quad G_{0}=0,
\end{equation}
in which the left-continuous volatility process $(V_{t})_{t\geq0}$ is defined by 
\begin{equation}\label{equ2.6}
V_{t}=\alpha_{0}+\mathbf{a}^{\prime}{\mathbf{Y}_{t-}},\qquad t>0,\qquad V_{0}=\alpha_{0}+\mathbf{a}^{\prime}{\mathbf{Y}_{0}},
\end{equation}
where the state process $(\mathbf{Y}_{t})_{t\geq0}$ is the unique c\`adl\`ag solution of the stochastic differential equation
\begin{equation}\label{equ2.7}
d\mathbf{Y}_{t}=B\mathbf{Y}_{t-}dt+\mathbf{e}(\alpha_{0}+\mathbf{a}^{\prime}\textbf{Y}_{t-})d[S,S]_{t},\qquad t>0,
\end{equation}
$d$ denotes differentiation with respect to $t$. The initial value $\mathbf{Y}_{0}$ is $\mathcal{F}_{0}$-measurable and independent of the driving semi-L\'evy process $(S_{t})_{t\geq0}$, and
\begin{align}\label{equ2.8}
B = \begin{pmatrix}
0 & 1 & 0 & \cdots & 0\\
0 & 0 & 1 & \cdots & 0\\
\vdots & \vdots & \vdots & \ddots & \vdots \\
0 & 0 & 0 & \cdots & 1\\
-\beta_{q} & -\beta_{q-1} & -\beta_{q-2} & \cdots & -\beta_{1}\\
\end{pmatrix},\quad
\mathbf{a}=\begin{pmatrix}
\alpha_{1}\\
\alpha_{2}\\
\vdots\\
\alpha_{q}\\
\end{pmatrix},\quad 
\mathbf{e}=\begin{pmatrix} 
0\\
0\\
\vdots\\
0\\
1\\
\end{pmatrix}.
\end{align}

\section{Periodic stationarity conditions}
\setcounter{equation}{0}
In this section we provide some conditions to prove that the volatility process $(V_{t})_{t\geq0}$ defined by 
(\ref{equ2.6}) is strictly periodically stationary with period $\tau$. As a result of main theorem, we prove that the
increments with constant length of process $(G_{t})_{t\geq0}$ is itself a periodically correlated (PC) process which
is the mian aim of this paper. We also give a sufficient an necessary condition by which we can determine the
volatility is non-negative. In the following theorem in ($b$) a $L^{r}-$matrix norm of the ($q\times q$)-matrix $C$ is
defined as 
\begin{align*}
\Vert C\Vert_{r}=\sup_{\mathbf{c}\in\mathbb{C}^{q}\setminus\lbrace0\rbrace}\frac{\Vert C\mathbf{c}\Vert_{r}}{\Vert
\mathbf{c}\Vert_{r}}.
\end{align*}
\begin{theorem}
(a) Let $(\mathbf{Y}_{t})_{t\geq0}$ be the state process of the SS-COGARCH(p,q) process with parameters $B$, 
$\mathbf{a}$ and $\alpha_{0}$ defined by (\ref{equ2.5}). Suppose that $(S_{t})_{t\geq0}$ be a simple semi-L\'evy
process defined by (\ref{equ2.2}). Then for all $0\leq s\leq t$ 
\begin{align}\label{equ3.1}
\mathbf{Y}_{t}=\mathfrak{J}_{s,t}\mathbf{Y}_{s}+\mathfrak{K}_{s,t},
\end{align}
where $(\mathfrak{J}_{s,t},\mathfrak{K}_{s,t})_{0\leq s\leq t}$ is a family of random $(q\times q)-$matrix
$\mathfrak{J}_{s,t}$ and random vector $\mathfrak{K}_{s,t}$ in $\mathbb{R}^{q}.$ In addition, 
$\big(\mathfrak{J}_{s+k\tau,t+k\tau},
\mathfrak{K}_{s+k\tau,t+k\tau}\big)_{k\in\mathbb{N}^{0}}$ are independent and identically distributed.\\

(b) Let $\eta_{i}, i=1,\cdots,q,$ be the eigenvalues of invertible matrix $B$ which have strictly negative real parts.
Also suppose that exists one $r\in[1,\infty]$ such that
\begin{align}\label{equ3.2}
\int_{\mathbb{R}}log\big(1+||P^{-1}\mathbf{e}\mathbf{a}^{\prime}P||_{r}z^{2}\big)d\nu_{t}(z)<-\frac{1}
{\Lambda(t+\tau)-\Lambda(t)}\nu_{t}(\mathbb{R})\eta\tau,\qquad \forall t\in[0,\tau),
\end{align}
where $P$ is a matrix in which $P^{-1}BP$ is diagonal and $\eta:=\max_{i=1,\cdots,q}\eta_{i}$ and 
$(\nu_{t})_{t\geq0}$ is semi-L\'evy measure defined by (\ref{equ2.3}). Then $\mathbf{Y}_{t+m\tau}$ converges in
distribution to a finite random vector $\mathbf{U}^{(t)}$ for fixed $t\in[0,\tau),$ as $m$ goes to infinity. The
distribution of the vector $\mathbf{U}^{(t)}$ is the unique solution of the random equation 
\begin{align}\label{equ3.3}
\mathbf{U}^{(t)}\overset{d}{=}\mathfrak{J}_{t,t+\tau}\mathbf{U}^{(t)}+\mathfrak{K}_{t,t+\tau},
\end{align}
where $\mathbf{U}^{(t)}$ is independent of $\big(\mathfrak{J}_{t,t+\tau},\mathfrak{K}_{t,t+\tau}\big)$.\\

(c) Let the conditions of (b) hold and $\mathbf{Y}_{0}\overset{d}{=}\mathbf{U}^{(0)}$, Then 
$(\mathbf{Y}_{t})_{t\geq0}$ and $(V_{t})_{t\geq0}$ are strictly periodically stationary with period $\tau$. In the
other hands, for any $s_{1}, s_{2}, \cdots, s_{n}\geq0$ and Borel sets $\mathbf{E}_{1}, \mathbf{E}_{2},\cdots, 
\mathbf{E}_{n}$ of $\mathbb{R}^{d}$ and Borel sets $J_{1}, J_{2}, \cdots, J_{n}$ of $\mathbb{R}$ and $k\in\mathbb{N}$,
\begin{align*}
P\big(\mathbf{Y}_{s_{1}}\in \mathbf{E}_{1}, \mathbf{Y}_{s_{2}}\in \mathbf{E}_{2}, \cdots, \mathbf{Y}_{s_{n}}\in
\mathbf{E}_{n}\big)=P\big(\mathbf{Y}_{s_{1}+k\tau}\in \mathbf{E}_{1}, \mathbf{Y}_{s_{2}+k\tau}\in \mathbf{E}_{2}, 
\cdots, \mathbf{Y}_{s_{n}+k\tau}\in \mathbf{E}_{n}\big),
\end{align*}
and
\begin{align*}
P\big(V_{s_{1}}\in J_{1}, V_{s_{2}}\in J_{2}, \cdots, V_{s_{n}}\in J_{n}\big)=P\big(V_{s_{1}+k\tau}\in J_{1},
V_{s_{2}+k\tau}\in J_{2}, \cdots, V_{s_{n}+k\tau}\in J_{n}\big).
\end{align*}
\end{theorem}
Proof: see Appendix P3.\\ \\
In the following remark we describe the non-negativity of the Lyapunov exponent which leads to the absolutely
convergence of the state process $(\mathbf{Y}_{t})_{t\geq0}$ in Theorem 3.1.
\begin{remark}
($a$) The proof of Theorem 3.1 will be based on the use of the general theory of multivariate random recurrence
equations, as discussed by Bougerol and Picard \cite{b6}, Brandt \cite{b7} and Vervaat \cite{v} in the one
dimensional
case. The state vector $(\mathbf{Y}_{t})_{t\geq0}$ defined by (\ref{equ2.7}) satisfies multivariate random recurrence
equation.\\

($b$) The condition (\ref{equ3.2}) which provides the stability of the model based on the existence of a vector norm
$||\cdot||_{r}$  such that $\mathfrak{J}_{t,t+\tau}$ and $\mathfrak{K}_{t,t+\tau}$ for all $t\in[0,\tau)$ satisfy the
conditions
\begin{align}\label{equ3.4}
E\Big(log||\mathfrak{J}_{t,t+\tau}||_{r}\Big)<0,\qquad E\Big(log^{+}||\mathfrak{K}_{t,t+\tau}||_{r}\Big)<\infty
\end{align}
where $log^{+}(x)=log(max\lbrace1,x\rbrace).$ $E\big(log||\mathfrak{J}_{t,t+\tau}||_{r}\big)<0$ is equivalent to the
assertion that the Lyapunov exponent of the $\big(\mathfrak{J}_{t+k\tau,t+(k+1)\tau}\big)_{k\in\mathbb{N}^{0}}$ is
strictly negative almost surely. i.e.
\begin{align*}
\limsup_{k\longrightarrow\infty}\frac{1}{k}log||\mathfrak{J}_{t,t+\tau}\cdots\mathfrak{J}_{t+k\tau,t+(k+1)\tau}||
_{r}<0,\quad a.s.
\end{align*}

($c$) The conditions of Theorem 3.1 imply (\ref{equ3.4}) with the natural matrix norm $||A||_{B,r}=||P^{-1}AP||_{r},$
for some matrix $A,$ which corresponds to the following the natural vector norm
\begin{align*}
||\mathbf{c}||_{B,r}:=||P^{-1}\mathbf{c}||_{r},\qquad \mathbf{c}\in\mathbb{C}^{q}
\end{align*}
where $P$ is a matrix in which $P^{-1}AP$ is diagonal.
\end{remark}
\begin{corollary}
If $(V_{t})_{t\geq0}$ is a strictly periodically stationary process with period $\tau$, then increments with constant
length of the process $(G_{t})_{t\geq0}$  make a PC process. In the other words, for any $t\geq0$ and $h\geq p>0$ and
$k\in\mathbb{N}$,
\begin{align*}
E\big(G_{t}^{(p)}\big)=E\big(G_{t+k\tau}^{(p)}\big),
\end{align*}
\begin{align*}
cov\big(G_{t}^{(p)},G_{t+h}^{(p)}\big)=cov\big(G_{t+k\tau}^{(p)},G_{t+h+k\tau}^{(p)}\big),
\end{align*}
where $G_{t}^{(p)}:=\int_{t}^{t+p}\sqrt{V_{s}}dS_{s}.$
\end{corollary}
Proof: see Appendix P4.
\begin{theorem}
Let $(\mathbf{Y}_{t})_{t\geq0}$ be the state process of the SS-COGARCH(p,q) process $(G_{t})_{t\geq0}$ with
parameters $B$, $\mathbf{a}$ and $\alpha_{0}>0.$ Suppose that $\gamma\geq-\alpha_{0}$ is a real constant and the
following two conditions hold:
\begin{align}\label{equ3.5}
\mathbf{a}^{\prime}e^{Bt}\mathbf{e}\geq0\quad \forall t\geq0,
\end{align}
\begin{align}\label{equ3.6}
\mathbf{a}^{\prime}e^{Bt}\mathbf{Y}_{0}\geq\gamma \quad a.s.\quad \forall t\geq0.
\end{align}
Then, with probability one,
\begin{align*}
V_{t}\geq\alpha_{0}+\gamma\geq0\quad\forall t\geq0.
\end{align*}
Conversely, if either (\ref{equ3.6}) fails, or (\ref{equ3.6}) holds with $\gamma>-\alpha_{0}$ and (\ref{equ3.5})
fails, then there exists a simple semi-L\'evy process $(S_{t})_{t\geq0}$ and $t_{0}\geq0$ such that $P(V_{t_{0}} < 0)
>0.$
\end{theorem}
The proof of the non-negativity volatility process $(V_{t})_{t\geq0}$ is similar to the
proof of Theorem 5.1 in \cite{b9} for L\'evy process.

\section{Characterization of the state process}
\setcounter{equation}{0}
The aim of this section is to study expected value and covariance function of the state process $\lbrace Y_{t}: t\geq0\rbrace$ and volatility process $\lbrace V_{t}:t\geq0\rbrace$. First, we prove that by some sufficient conditions the expected value and covariance $\mathbf{Y}_{t}$ exist. Then, by presenting the first and second moments of the random vector $U^{(0)},$ we find the expected value and covariance function of the state process. Furthermore, a closed form for square increments of the COGARCH process is characterized. 
\begin{lemma}
Let the assumptios of Theorem 3.1 hold. If $E\big(||\mathbf{Y}_{0}||_{r}\big)^{c}<\infty,$  for $c=1,2,$ then

($a$) If $E(S_{t}^{2})<\infty,$ then $E(\mathbf{Y}_{t})<\infty$ and $E(\mathbf{U}^{(0)})<\infty.$

($b$)  If $E(S_{t}^{4})<\infty,$ then $cov(\mathbf{Y}_{t})<\infty$ and $cov(\mathbf{U}^{(0)})<\infty.$

where $\lbrace S_{t}: t\geq0\rbrace$ is the simple semi-Levy process.
\end{lemma}
Proof: see Appendix P5.
\begin{remark}
By Theorem 3.1(b), ($a$) we find that $E(\mathbf{U}^{(0)})$ is the solution of the following random equation
\begin{align*}
\big(I-E(\mathfrak{J}_{0,\tau})\big)E\big(\mathbf{U}^{(0)}\big)=E(\mathfrak{K}_{0,\tau}),
\end{align*}
and

($b$) $E\big((\mathbf{U}^{(0)})(\mathbf{U}^{(0)})^{\prime}\big)$ is the solution of the following equation
\begin{align*}
\big(I_{q^{2}}-E(\mathfrak{J}_{0,\tau}\otimes\mathfrak{J}_{0,\tau})\big)vec\Big(E\big((\mathbf{U^{(0)}})(\mathbf{U}^{(0)})^{\prime}\big)\Big)
&=\big(E(\mathfrak{K}_{0,\tau}\otimes\mathfrak{J}_{0,\tau})+E(\mathfrak{J}_{0,\tau}\otimes\mathfrak{K}_{0,\tau})\big)E(\mathbf{U}^{(0)})\\
&\quad+vec\big(E(\mathfrak{K}_{0,\tau}\mathfrak{K}_{0,\tau}^{\prime})\big),
\end{align*}

where $\otimes$ is the Kronecker product of two matrices and for a matrix $\mathbf{C},$ $vec(\mathbf{C})$ is the column vector in
$\mathbb{C}^{q^{2}}$ which is constructed by stacking the columns of matrix $A$ in a vector.
\end{remark}
The following lemmas establish the the mean and covariance function of the state process. 
\begin{lemma}
Suppose that $\lbrace \mathbf{Y}_{t}:t\geq0\rbrace$ be the state process and the conditions of Theorem 3.1 and lemma 4.1 hold. Then for  $t, h\geq0,$ there exists $m, n\in\mathbb{N}^{0}$ and $t_{1}, t_{2}\in[0,\tau)$ such that $t\in\big[m\tau,(m+1)\tau\big),$ $t+h\in\big[n\tau,(n+1)\tau\big),$ $t=t_{1}+m\tau$ and $t+h=t_{2}+n\tau$ and
 \begin{align}\label{equ4.1}
E(\mathbf{Y}_{t})=E(\mathfrak{J}_{0,t_{1}})E(\mathbf{U})+E(\mathfrak{K}_{0,t_{1}}),
\end{align}
\begin{align}\nonumber
cov(\mathbf{Y}_{t},\mathbf{Y}_{t+h})=E(\mathfrak{J}_{0, t_{2}})\big(E(\mathfrak{J}_{0,\tau})\big)^{m-n-1}
&\Big[E(\mathfrak{J}_{0,\tau}E(\mathbf{U}\mathbf{U}^{\prime})\mathfrak{J}^{\prime}_{0,t_{1}})-E(\mathfrak{J}_{0,\tau})E(\mathbf{U})E(\mathbf{U}^{\prime})E(\mathfrak{J}^{\prime}_{0,t_{1}})\\ \nonumber
&+E(\mathfrak{J}_{0,\tau}E(\mathbf{U})\mathfrak{K}^{\prime}_{0,t_{1}})-E(\mathfrak{J}_{0,\tau})E(\mathbf{U})E(\mathfrak{K}^{\prime}_{0,t_{1}})\\ \nonumber
&+E(\mathfrak{K}_{0,\tau}E(\mathbf{U^{\prime}})\mathfrak{J}^{\prime}_{0,t_{1}})-E(\mathfrak{K}_{0,\tau})E(\mathbf{U^{^{\prime}}})E(\mathfrak{J}^{\prime}_{0,t_{1}})\\ \label{equ4.2}
&+E(\mathfrak{K}_{0,\tau}\mathfrak{K}^{\prime}_{0,t_{1}})-E(\mathfrak{K}_{0,\tau})E(\mathfrak{K}^{\prime}_{0,t_{1}})\Big].
\end{align}
\end{lemma}
Proof: see Appendix P6.
\begin{corollary}
Let $\lbrace V_{t}: t\geq0\rbrace$ be the volatility process. Then for $t, h\geq0,$ expected value and covariance function of $V_{t}$ have the following forms.
$$E(V_{t})=\alpha_{0}+\mathbf{a}^{\prime}E(\mathbf{Y}_{t})$$
$$cov(V_{t}, V_{t+h})=\mathbf{a}^{\prime}cov(\mathbf{Y}_{t},\mathbf{Y}_{t+h})\mathbf{a}.$$
\end{corollary}
Proof: see Appendix P7.\\

In financial time series, the returns have negligible correlation while the squared returns are significantly correlated, therefore we investigate the behavior of the second-order properties of the increments of the COGARCH process. We assume that volatility process  is strictly periodically stationary and non-negative.\\
Now we present the first and second orders of the increment process $G_{t}^{(p)}$ that in defined in Corollary 4.2.

\begin{proposition}
Let $G$ be a zero mean simple semi-Levy driven COGARCH process. Then for $t\geq0$ and $h\geq p>0,$

($a$)
\begin{align}\label{equ4.3}
E(G_{t}^{(p)}\big)&=0,\\ \label{equ4.4}
cov\big(G_{t}^{(p)},G_{t+h}^{(p)}&\big)=0.
\end{align}

($b$) There exist $m, m^{\prime}\in\mathbb{N}^{0}$ where $t\in A_{ml+i}, i=1, \cdots, l$ and $t+p\in A_{(m+m^{\prime})l+i^{\prime}}, i^{\prime}=i,\cdots, l,$ then
\begin{align}\nonumber
E\big((G_{t}^{(p)})^{2}\big)&=E\big((Z_{i})^{2}\big)\int_{t}^{t_{ml+i}}E\big(V_{s}\big)\lambda(s)ds+E\big((Z_{i^{\prime}})^{2}\big)
\int_{t_{(m+m^{\prime})l+i^{\prime}-1}}^{t+p}E\big(V_{s}\big)\lambda(s)ds\\ \label{equ4.5}
&\quad+\sum_{r=i}^{m^{\prime}l+i^{\prime}-2}E\big((Z_{r+1})^{2}\big)\int_{t_{ml+r}}^{t_{ml+r+1}}E\big(V_{s}\big)\lambda(s)ds.
\end{align}
Moreover, there exist $n, n^{\prime}\in\mathbb{N}^{0}$ ($n\geq m$ and $n^{\prime}\geq m^{\prime}$) where $t+h\in A_{nl+j}, j=1, \cdots, l$ and $t+h+p\in A_{(n+n^{\prime})l+j^{\prime}}, j^{\prime}=j,\cdots, l,$ then
\begin{align}\nonumber
cov\big((G_{t}^{(p)})^{2},(G_{t+h}^{(p)})^{2}\big)&=
E\big((Z_{j})^{2}\big)\mathbf{a}^{\prime}\int_{t+h}^{t_{nl+j}}\lambda(s)E(\mathfrak{J}_{t+p,s^{-}})cov\big((G_{t}^{(p)})^{2},\mathbf{Y}_{t+p}\big)ds\\ \nonumber
&\quad+E\big((Z_{j^{\prime}})^{2}\big)\mathbf{ a}^{\prime}
\int_{t_{(n+n^{\prime})l+j^{\prime}-1}}^{t+h+p}\lambda(s)E(\mathfrak{J}_{t+p,s^{-}})cov\big((G_{t}^{(p)})^{2},\mathbf{Y}_{t+p}\big)ds\\ \label{equ4.6}
&\quad+\sum_{r=j}^{n^{\prime}l+j^{\prime}-2}E\big((Z_{r+1})^{2}\big)\mathbf{a}^{\prime}\int_{t_{nl+r}}^{t_{nl+r+1}}\lambda(s)E(\mathfrak{J}_{t+p,s^{-}})
cov\big((G_{t}^{(p)})^{2},\mathbf{Y}_{t+p}\big)ds.
\end{align}
\end{proposition}
Proof: see Appendix P8.
\begin{remark}
For $s\geq0,$ if we assume that $s\geq0,$ $\int_{\mathbb{R}}z^{3}\nu_{s}(dz)=0,$ Then
\begin{align*}
cov\big((G_{t}^{(p)})^{2},\mathbf{Y}_{t+p}\big)&=
2E(I_{t+p}\mathbf{Y}_{t+p})-2E(\mathfrak{J}_{t,t+p})E(I_{t}\mathbf{Y}_{t})\\
&\quad-\big(cov(\mathbf{Y}_{t+p})+cov(\mathbf{Y}_{t+p},\mathbf{Y}_{t})-\int_{t^{+}}^{t+p}cov(\mathbf{Y_{t+p}},\mathbf{Y}_{s^{-}})dsB^{\prime}\big)\mathbf{e},
\end{align*}  
where $I_{t}:=\int_{0}^{t}\big(G_{s^{-}}\sqrt{V_{s}}\big)dS_{s}$ and
$$E(I_{t}\mathbf{Y}_{t})=B\int_{0}^{t}E(I_{s}\mathbf{Y}_{s})ds+\alpha_{0}\int_{0}^{t}\int_{\mathbb{R}}
E(I_{s}\mathbf{Y}_{s})z^{2}\nu_{s}(dz)ds.$$
\end{remark}

\section{Simulation}
\setcounter{equation}{0}
In this section we simulate the simple semi-L\'evy process defined by (\ref{equ2.2}). This process is a compound Poisson process with time-varying arrival rate $\Lambda(t)$ defined by (\ref{equ2.1}). Then we verify the theoretical results concerning PC structure of the increments of the SS-COGARCH(p,q) process $(G_{t})_{t\geq0}$ defined by (\ref{equ2.5}) by simulation. For this we simulate the state process $(\mathbf{Y}_{t})_{t\geq0}$ defined by (\ref{equ2.7}) at jump time points and non-jump time points using its random recurrence equation (\ref{equ3.1}). Then we evaluate the discretized version of the volatility process $(V_{t})_{t\geq0}$ defined by (\ref{equ2.6}) and corresponding the SS-COGARCH(p,q) process $(G_{t})_{t\geq0}$. Finally, we verify the PC structure of the increments of the SS-COGARCH process by following the method of \cite{d2}.\\
For simulating the simple semi-L\'evy process defined by (\ref{equ2.2}) with the underlying Poisson process $\big(N(t)\big)_{t\geq0},$ we consider $T_{1}$ as the time of the first jump and $T_{n}, n=2, 3, \cdots$ the time intervals between the $(n-1)^{th}$ and $n^{th}$  jumps. Then $\Upsilon_{n}=\sum_{j=1}^{n}T_{j}, n\in\mathbb{N},$ are the arrival times and $\Upsilon_{0}=0$. Therefore for $j=1, 2, \cdots$
\begin{align}\nonumber
F_{T_{n}}^{s}(x)&:=P\big(T_{n}\leq x\big{|}\Upsilon_{n-1}=s\big)\\ \nonumber
&=1-P\big(N(s+x)-N(s)=0\big)\\\label{equ5.1}
&=1-e^{-\Lambda(s+x)+\Lambda(s)},
\end{align}
where $\Lambda(t)$ defined by (\ref{equ2.1}). The arrival times $\Upsilon_{1}, \Upsilon_{2}, \cdots$ are generated by the following algorithm.

1. Generate the independent and identically distributed (iid) sequence $U_{1}, U_{2}, \cdots$ from Uniform (0,1). Then by (\ref{equ5.1}) as $\Upsilon_{0}=0$ the first arrival time $\Upsilon_{1}=T_{1}$ has distribution $F_{\Upsilon_{1}}^{0}(x)=1-e^{-\Lambda(x)}.$ Therefore
$$\Lambda(\Upsilon_{1})\overset{d}{=}-ln(1-U),$$
where $U$ denotes a Uniform (0,1). So by generating $U_{1}$, $\Upsilon_{1}=\Lambda^{-1}\big(-ln(1-U_{1})\big)$ can be considered as a generated sample for the first arrival time. If $\Upsilon_{n-1}, n=2, 3, \cdots,$ is the $(n-1)^{th}$ evaluated arrival time, then by (\ref{equ5.1}) $T_{n}$ has distribution $F_{T_{n}}^{\Upsilon_{n-1}}(x)=1-e^{-\Lambda(x+\Upsilon_{n-1})+\Lambda(\Upsilon_{n-1})}.$ Therefore
$$\Lambda(\Upsilon_{n})\overset{d}{=}\Lambda(\Upsilon_{n-1})-ln(1-U).$$
So by generating $U_{n},$ $\Upsilon_{n}=\Lambda^{-1}\big(\Lambda(\Upsilon_{n-1})-ln(1-U_{n})\big)$ is a generated sample for the $n^{th}$
arravial time. Thus applying the iid sample $U_{1}, U_{2}, \cdots$ we can evaluate successively the $n^{th}$ arrival time by the (\ref{equ5.1}),
for the details see \cite[p.99]{c1}. So by having the periodic intensity function $\lambda(u)$ in (\ref{equ2.1}), one can evaluate $\Lambda^{-1}
(\cdot)$ by available software.

2. Consider some periodic drift function $H_{t}$ and as the successive jump size $Z_{n}$ generate independently and has distribution $F_{n}
(\cdot)$ if corresponding arrival time belongs to $\mathfrak{D}_{j}=\bigcup_{k=0}^{\infty}A_{j+kl}, j=1, 2, \cdots,l.$ Now evaluate the simple
semi-L\'evy process $(S_{t})_{t\geq0}$ from (\ref{equ2.2}) as
\begin{align}\label{equ5.2}
S_{t}=H_{t}+\sum_{n=1}^{N(t)}\sum_{j=1}^{l}Z_{n}^{j}I_{\lbrace\Upsilon_{n}\in\mathfrak{D}_{j}\rbrace}.
\end{align}

Now we consider the following steps for the simulation of the SS-COGARH(p,q) process defined by (\ref{equ2.5})-(\ref{equ2.7}).

1. Consider $p$ and $q$ as some integer such that $q\geq p\geq1.$

2. Choose real parameters $\beta_{1}, \cdots, \beta_{q}$ and $\alpha_{1}, \cdots, \alpha_{p}$ and $\alpha_{0}>0$ such that the eigenvalues of the matrix $B$ defined by (\ref{equ2.8}) have strictly negative real parts and conditions (\ref{equ3.2}), (\ref{equ3.5}) and (\ref{equ3.6}) are satisfied.

3. Having evaluated arrival times $\Upsilon_{n}$ by the above algorithm, generate the state process $(\mathbf{Y}_{\Upsilon_{n}})_{n\in\mathbb{N}}$ by the following the recurrence equation after assuming some initial value for $\mathbf{Y}_{\Upsilon_{0}}$
\begin{align*}
\mathbf{Y}_{\Upsilon_{n}}=e^{B\big(\Upsilon_{n}-\Upsilon_{n-1}\big)}\mathbf{Y}_{\Upsilon_{n-1}}+\mathbf{e}\Big(\alpha_{0}+
\mathbf{a}^{\prime}e^{B\big(\Upsilon_{n}-\Upsilon_{n-1}\big)}\mathbf{Y}_{\Upsilon_{n-1}}\Big)\big(Z_{n}\big)^{2}, \quad n\in\mathbb{N}.
\end{align*}
This recurrence equation obtained by replacing $s=\Upsilon_{n-1}$ and $t=\Upsilon_{n}$ in (\ref{equ3.1}). The jump size $Z_{n}$ can be
simulated by (\ref{equ5.2}) for predefined distributions $\mu_{1}, \mu_{2}, \cdots,\mu_{l}.$

4. As the simple semi-L\'evy process $(S_{t})_{t\geq0}$ (\ref{equ2.2}) has no jump over $[\Upsilon_{n-1},\Upsilon_{n}^{-}], n\in\mathbb{N},$
it follows from (\ref{equ2.7}) that $d\mathbf{Y}_{t}=B\mathbf{Y}_{t}dt$ for $t\in[\Upsilon_{n-1},\Upsilon_{n}^{-}].$ Therefore for
$t\in[\Upsilon_{n-1},\Upsilon_{n}^{-}]$
\begin{align*}
e^{-Bt}d\mathbf{Y}_{t}=Be^{-Bt}\mathbf{Y}_{t}dt
\end{align*}
so that
\begin{align*}
d\big(e^{-Bt}\mathbf{Y}_{t}\big)=0.
\end{align*}
From this follows that for $t\in[\Upsilon_{n-1},\Upsilon_{n}^{-}]$
\begin{align*}
\int_{\Upsilon_{n-1}}^{t}d\big(e^{-Bu}\mathbf{Y}_{u}\big)=0
\end{align*}
hence
\begin{align}\label{equ5.3}
\mathbf{Y}_{t}=e^{B(t-\Upsilon_{n-1})}\mathbf{Y}_{\Upsilon_{n-1}}.
\end{align}
By (\ref{equ2.6}) and (\ref{equ5.3}), the discrete-time version of the process $(V_{t})_{t\geq0}$ is as
\begin{align}\nonumber
V_{\Upsilon_{n}}&=\alpha_{0}+\mathbf{a}^{\prime}\mathbf{Y}_{\Upsilon_{n}^{-}}\\ \label{equ5.4}
&=\alpha_{0}+\mathbf{a}^{\prime}e^{B(\Upsilon_{n}-\Upsilon_{n-1})}\mathbf{Y}_{\Upsilon_{n-1}},
\end{align}
and using that the process $(S_{t})_{t\geq0}$ (\ref{equ2.2}) has one jump at time $\Upsilon_{n}$ over $[\Upsilon_{n-1},\Upsilon_{n}]$ it follows from (\ref{equ2.5}) that 
\begin{align}\nonumber
G_{\Upsilon_{n}}-G_{\Upsilon_{n-1}}&=\int_{0}^{\Upsilon_{n}}\sqrt{V_{u}}dS_{u}-\int_{0}^{\Upsilon_{n-1}}\sqrt{V_{u}}dS_{u}\\ \nonumber
&=\int_{\Upsilon_{n-1}}^{\Upsilon_{n}}\sqrt{V_{u}}dS_{u}\\ \label{equ5.5}
&=\sqrt{V_{\Upsilon_{n}}}Z_{n}.
\end{align}

5. Having evaluated values of the process $(\mathbf{Y}_{\Upsilon_{n}})_{n\in\mathbb{N}}$ and $G_{0}=0,$ generate the process $(V_{\Upsilon_{n}})_{n\in\mathbb{N}}$ by (\ref{equ5.4}) and corresponding the process $(G_{\Upsilon_{n}})_{n\in\mathbb{N}}$ by (\ref{equ5.5}).

6. Finally, using the values of $V_{\Upsilon_{n}}$ and $G_{\Upsilon_{n}}$ provided by previous step, evaluate the sampled processes $(V_{ih})_{i\in\mathbb{N}}$ and $(G_{ih})_{i\in\mathbb{N}}$ for some  $h>0$ by the followings:

(i) Suppose that  $ih\in[\Upsilon_{n-1},\Upsilon_{n}),$ for $i, n\in\mathbb{N}.$ Since the simple semi-L\'evy process $(S_{t})_{t\geq0}$ (\ref{equ2.2}) has no jump over $[\Upsilon_{n-1},\Upsilon_{n}),$ it follows from (\ref{equ2.6}) and (\ref{equ5.3}) that for $ih\in[\Upsilon_{n-1},\Upsilon_{n})$
\begin{align*}
V_{ih}&=\alpha_{0}+\mathbf{a}^{\prime}\mathbf{Y}_{ih^{-}}\\
&=\alpha_{0}+\mathbf{a}^{\prime}e^{B(ih-\Upsilon_{n-1})}\mathbf{Y}_{\Upsilon_{n-1}},
\end{align*}
note that if $ih=\Upsilon_{n-1},$ then it follows from step 4 that
\begin{align*}
V_{ih}&=\alpha_{0}+\mathbf{a}^{\prime}\mathbf{Y}_{\Upsilon_{n-1}^{-}}\\
&=\alpha_{0}+\mathbf{a}^{\prime}e^{B(\Upsilon_{n-1}-\Upsilon_{n-2})}\mathbf{Y}_{\Upsilon_{n-2}}.
\end{align*}

(ii) Using that the process $(S_{t})_{t\geq0}$ (\ref{equ2.2}) has no jump over $[\Upsilon_{n-1},ih]$ it follows from (\ref{equ2.5}) that
\begin{align*}
G_{ih}-G_{\Upsilon_{n-1}}&=\int_{0}^{ih}\sqrt{V_{u}}dS_{u}-\int_{0}^{\Upsilon_{n-1}}\sqrt{V_{u}}dS_{u}\\
&=\int_{\Upsilon_{n-1}}^{ih}\sqrt{V_{u}}dS_{u}=0\\
\end{align*}
hence
\begin{align*}
G_{ih}=G_{\Upsilon_{n-1}}.
\end{align*}

\subsection{Test for the PC Structure of the increments process}
To detect the PC structure of a process, Hurd and Miamee \cite{h1} and Dudek et al. \cite{d2} showed that their proposed spectral coherence can be used to test whether a discrete-time process is PC. Their method is based on the fact that the support of the spectral coherence of a PC process with period $\varrho$ is contained in the subset of parallel lines $\lambda_{s}=\lambda_{r}+2j\pi/\varrho$ for $j=-(\varrho-1), \cdots, -1, 0, 1, \cdots, (\varrho-1).$ The squared coherence statistic for the series $X_{1}, X_{2}, \cdots, X_{N}$ is computed as follows
\begin{align*}
|\hat{\gamma}(\lambda_{r},\lambda_{s},M)|^{2}=\frac{|\sum_{m=1}^{M-1}\tilde{X}_{N}(\lambda_{r-M/2+m})\overline{\tilde{X}_{N}(\lambda_{s-M/2+m})}|^{2}}{\sum_{m=1}^{M-1}|\tilde{X}_{N}(\lambda_{r-M/2+m})|^{2}\sum_{m=1}^{M-1}|\tilde{X}_{N}(\lambda_{s-M/2+m})|^{2}}
\end{align*}
where $\tilde{X}_{N}(\lambda_{j})=\sum_{k=1}^{N}X_{k}e^{-i\lambda_{j}k}$ is discrete Fourier transform of $X_{k}$ for $j=0, 1, \cdots, N-1,$ $\lambda_{j}=2\pi j/N$ and $\lambda_{j}\in(0,2\pi].$  This statistic satisfies $0\leq|\hat{\gamma}(\lambda_{r},\lambda_{s},M)|^{2}\leq1.$ \\
Under the null hypothesis that $\tilde{X}_{N}(\lambda_{j})$ are complex Gaussian with uncorrelated real and imaginary parts for each $j$,  squared coherence statistic has probability density, \cite{h1}
\begin{align*}
p(|\gamma |^{2})=(M-1)\big(1-|\gamma|^{2}\big)^{M-2},\quad 0\leq |\gamma|^{2}\leq1.
\end{align*}
For type I error $\alpha,$ the squared coherence $\alpha$-threshold is determined, \cite{h1}
\begin{align*}
x_{\alpha}:=|\gamma|_{\alpha}^{2}=1-e^{log(\alpha)/(M-1)}.
\end{align*} 
The values of statistic $|\hat{\gamma}(\lambda_{r},\lambda_{s},M)|^{2}$ are computed for all $r$ and $s$ that pair $(\lambda_{r}, \lambda_{s})\in[0,2\pi)\times[0,2\pi).$ By plotting the values of statistic that exceed the $\alpha-$threshold, if there are some significant values of statistic that lie along the parallel equally spaced diagonal lines, then $X_{k}$ is PC. The graph of these significant values indicates the presence of the subset of parallel lines $s=r+jN/\varrho$ for $j=-(\varrho-1), \cdots, -1, 0, 1, \cdots, (\varrho-1).$ \\
To ensure that periodic structure of the series $X_{k}$ is not a consequence of a periodic mean, it is recommended to remove the periodic mean from this series first.
\begin{example}
Let $(S_{t})_{t\geq0}$ be a simple semi-L\'evy process by the rate function $\lambda(t)=-cos(\frac{\pi}{6}t)$ $+4,$ for $t\geq0.$ Furthermore,
$\tau=12,$ $l=5$ and the length of the successive partitions of each period intervals are 2, 2, 2, 3, 3. Moreover, the distribution of jumps size on
these subintervals are assumed to be $N(3,1),$ $N(0,1),$ $N(1.25,1.25),$ $N(4,1),$ and $N(0,1.5),$ where $N(\mu,\sigma^{2})$ denotes a
Normal distribution with mean $\mu$ and variance $\sigma^{2}$.
\end{example}
In this example we consider SS-COGARCH(1,3) process with parameters of $\alpha_{0}=1,$ $\alpha_{1}=0.03,$ $\beta_{1}=5,$ $\beta_{2}=9$
and $\beta_{3}=5.$  Thus, the matrix $B$ is 
$$
B = \begin{pmatrix}
0 & 1 & 0 \\
0 & 0 & 1 \\
-5 & -9 & -5\\
\end{pmatrix}
$$
and conditions (\ref{equ3.2}), (\ref{equ3.5}) and (\ref{equ3.6}) are satisfy. For such the SS-COGARCH process we simulate $G_{\Upsilon_{n}}$ for the duration of $40$ period intervals with the parameters specified above, $\mathbf{Y}_{0}=(8.3580, 2.3377, 0.9040)^{\prime}$ and $G_{0}=0$. Then, using step 6, we sample from this process in equally space partition with distance one (h=1). So we get 480 discretized samples of this 40 period intervals. Then we follow to verify that the increments of the sampled process are a PC process.
 \begin{figure}[H]
\centering
\includegraphics[scale=0.5]{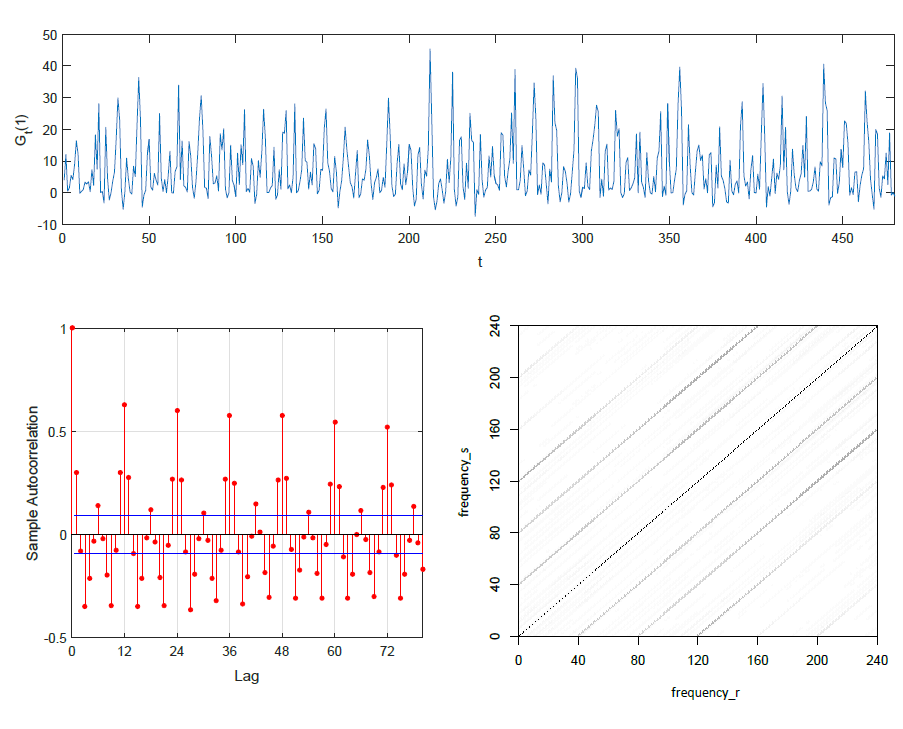}
\caption{\scriptsize{Top: the increments of the simulated process $\lbrace G_{i}:i\in\mathbb{N}\rbrace$ of size 480; bottom left: the sample autocorrelation plot of $G_{t}^{(1)}$; bottom right:  the significant values of the sample spectral coherence with $\alpha=0.05.$}}
 \end{figure}
In figure 1 graph of the increments of the sampled process of size 480 (top) with the sample autocorrelation graph of this process (bottom left) are presented. The bottom right graph shows the sample coherent statistics values for a specified collection of pairs $(\lambda_{r},\lambda_{s})\in[0,2\pi)\times[0,2\pi)$ and $M=240$ that exceed the threshold corresponding to $\alpha=0.05.$ The parallel lines for the sample spectral coherence confirm the increments of the sampled process are PC. Also in this graph, the significant off-diagonal is at $|r-s|=40$ which verifies the first peak at 40 and shows that there is a second order periodic structure with period $\varrho=480/40=12.$ 
\begin{table}[H]
\caption{\scriptsize{Some different values of $|\hat{\gamma}(\lambda_{r},\lambda_{s},M)|^{2}$ for some values $r, s$ with $\alpha=0.05$ and $x_{\alpha}=0.0125.$}}
\centering
\includegraphics[scale=0.4]{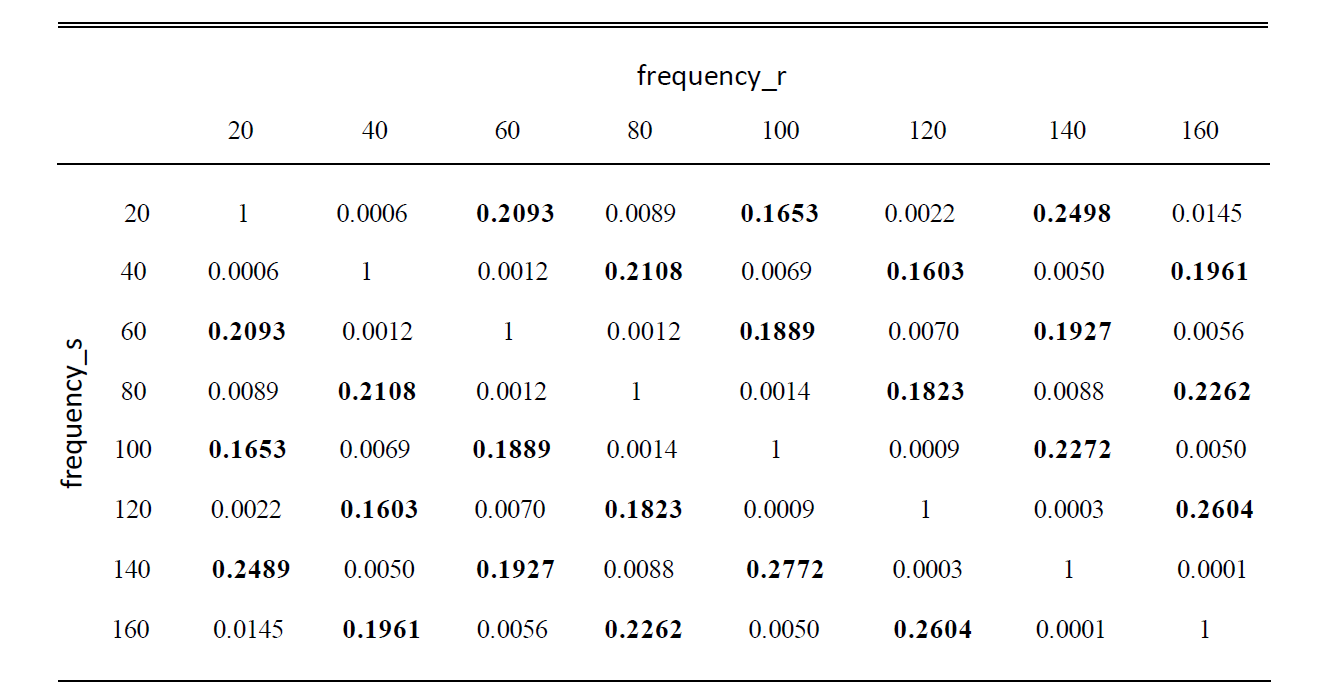}
\end{table}
The some different values of sample coherence statistics for the test that the increments of the sampled process have period 12 are presented in
Table 1. As the corresponding $\alpha=0.05$ threshold shows the test is significant on the corresponding parallel lines of Figure 1.

\section{Appendix}
\setcounter{equation}{0}
\textbf{P1: Proof of Lemma 2.2}\\
For any $t\geq0$ there exist $j=1,\cdots, l$ and $m\in\mathbb{N}^{0}$ such that $t\in A_{j+ml}$. Thus, using the 
Definition 2.2, Definition 2.3 and the fact that $(S_{t})_{t\geq0}$ has independence increments we have
\begin{align*}
E\big(e^{iwS_{t}}\big)&=E\big(e^{iw\big(D_{t}+\sum_{n=1}^{N(t)}Z_{n}\big)}\big)=e^{iwD_{t}}E\big(e^{iw\sum_{n=1}^{N(t)}Z_{n}}\big)\\
&=e^{iwD_{t}}E\Big(e^{iw\sum_{n=1}^{N(t_{1})}Z_{n}^{1}}\times 
e^{iw\sum_{n=N(t_{1})+1}^{N(t_{2})}Z_{n}^{2}}\times\cdots\times e^{iw\sum_{n=N(t_{ml+j-1})+1}^{N(t)}Z_{n}^{j}}\Big)\\
&=e^{iwD_{t}}\Big(\prod_{k=0}^{m-1}\prod_{r=1}^{l}E\big(e^{iw\sum_{n=N(t_{kl+r-1})+1}^{N(t_{kl+r})}Z_{n}^{r}}\big)\Big)
\times\Big(\prod_{r=1}^{j-1}E\big(e^{iw\sum_{n=N(t_{ml+r-1})+1}^{N(t_{ml+r})}Z_{n}^{r}}\big)\Big)\\
&\quad\times\Big(E\big(e^{iw\sum_{n=N(t_{ml+j-1})+1}^{N(t)}Z_{n}^{j}}\big)\Big).
\end{align*}
Since for $r=1,\cdots, l,$ $Z_{n}^{r}$ are independent and have distribution $F_{r},$ it follows from Definition 2.2 
and conditional expected value
for $k=0, \cdots, m$ and $r=1, \cdots, l$ that
\begin{align*}
E\big(e^{iw\sum_{n=N(t_{kl+r-1})+1}^{N(t_{kl+r})}Z_{n}^{r}}\big)&=E\Big(E\big(e^{iw\sum_{n=1}^{N(t_{kl+r})-
N(t_{kl+r-1})}Z_{n}^{r}}\vert N(t_{kl+r})-N(t_{kl+r-1})=N\big)\Big)\\
&=\sum_{N=0}^{\infty}E\big(e^{iw\sum_{n=1}^{N}Z_{n}^{r}}\big)P\big(N(t_{kl+r})-N(t_{kl+r-1})=N\big)\\
&=\sum_{N=0}^{\infty}\big(E(e^{iwZ_{n}^{r}})\big)^{N}\frac{\big(\Lambda(t_{kl+r})-\Lambda(t_{kl+r-1})\big)^{N}
e^{-\big(\Lambda(t_{kl+r})-\Lambda(t_{kl+r-1})\big)}}{N!}\\
&=e^{-\big(\Lambda(t_{kl+r})-\Lambda(t_{kl+r-1})\big)}\sum_{N=0}^{\infty}\frac{\Big[\int_{\mathbb{R}}e^{iwz}
\big(\Lambda(t_{kl+r})-\Lambda(t_{kl+r-1})\big)F_{r}(dz)
\Big]^{N}}{N!}\\
&=e^{\int_{\mathbb{R}}(e^{iwz}-1)\big(\Lambda(t_{kl+r})-\Lambda(t_{kl+r-1})\big)F_{r}(dz)}
\end{align*}
Therefore 
\begin{align*}
E\big(e^{iwS_{t}}\big)&=e^{iwD_{t}}e^{\int_{\mathbb{R}}
(e^{iwz}-1)\Big[\sum_{k=0}^{m-1}\sum_{r=1}^{l}\big(\Lambda(t_{kl+r})-\Lambda(t_{kl+r-1})\big)F_{r}(dz)}\\
&\qquad\qquad\qquad\qquad^{+\sum_{r=1}^{j-1}\big(\Lambda(t_{ml+r})-\Lambda(t_{ml+r-1})\big)F_{r}
(dz)+\big(\Lambda(t)-\Lambda(t_{ml+j-1})\big)F_{j}(dz)\Big]}\\
&=e^{iw\Big(D_{t}+\int_{|z|
\leq1}z\Big[\sum_{k=0}^{m-1}\sum_{r=1}^{l}\big(\Lambda(t_{kl+r})-\Lambda(t_{kl+r-1})\big)F_{r}(dz)}\\
&\qquad\qquad\quad\quad\quad^{+\sum_{r=1}^{j-1}\big(\Lambda(t_{ml+r})-\Lambda(t_{ml+r-1})\big)F_{r}
(dz)+\big(\Lambda(t)-\Lambda(t_{ml+j-1})\big)F_{j}(dz)\Big]\Big)}\\
&\quad\times e^{\int_{\mathbb{R}}(e^{iwz}-1-iwzI_{\lbrace|\mathit{z}|
\leq1\rbrace})\Big[\sum_{k=0}^{m-1}\sum_{r=1}^{l}\big(\Lambda(t_{kl+r})-\Lambda(t_{kl+r-1})\big)F_{r}(dz)}\\
&\qquad\qquad\qquad\qquad\quad\quad\quad^{+\sum_{r=1}^{j-1}\big(\Lambda(t_{ml+r})-\Lambda(t_{ml+r-1})\big)F_{r}
(dz)+\big(\Lambda(t)-\Lambda(t_{ml+j-1})\big)F_{j}(dz)\Big]}.
\end{align*}
\textbf{P2: Proof of Corollary 2.3}\\
It is sufficient to prove that for any $0\leq s<t$ and $K\in\mathbb{N}$,
$$S_{t}-S_{s}\overset{d}{=}S_{t+K\tau}-S_{s+K\tau}.$$
For any $0\leq s<t$ there exist $m, m^{\prime}\in\mathbb{N}^{0}$ and $j, j^{\prime}=1, \cdots, l$ such that $s\in 
A_{j+ml}$ and $t\in 
A_{j^{\prime}+(m+m^{\prime})l}.$ Thus
\begin{align*}
E\big(e^{iw\big(S_{t}-S_{s}\big)}\big)&=E\big(e^{iw\big(D_{t}-
D_{s}+\sum_{n=1}^{N(t)}Z_{n}-\sum_{n=1}^{N(s)}Z_{n}\big)}\big)
=e^{iw\big(D_{t}-D_{s}\big)}E\big(e^{iw\sum_{n=N(s)+1}^{N(t)}Z_{n}}\big)\\
&=e^{iw\big(D_{t}-D_{s}\big)}E\Big(e^{iw\sum_{n=N(s)+1}^{N(t_{ml+j})}Z_{n}^{j}}
\times\cdots\times e^{iw\sum_{n=N(t_{(m+m^{\prime})l+j^{\prime}-1})+1}^{N(t)}Z_{n}^{j^{\prime}}}\Big)\\
\end{align*}
\begin{align*}
&=e^{iw(D_{t}-D_{s})}\times E\big(e^{iw\sum_{N(s)+1}^{N(t_{ml+j})}Z_{n}^{j}}\big)
\times\Big(\prod_{r=j+1}^{l}E\big(e^{iw\sum_{n=N(t_{ml+r-1})+1}^{N(t_{ml+r})}Z_{n}^{r}}\big)\Big)\\
&\quad\times\Big(\prod_{k=m+1}^{m+m^{\prime}-1}\prod_{r=1}^{l}E\big(e^{iw\sum_{n=N(t_{kl+r-1})+1}^{N(t_{kl+r})}Z_{n}^{r}}\big)\Big)\\
&\quad\times\Big(\prod_{r=1}^{j^{\prime}-1}E\big(e^{iw\sum_{n=N(t_{(m+m^{\prime})l+r-1})+1}^{N(t_{(m+m^{\prime})l+r})}Z_{n}^{r}}\big)\Big)
\times E\big(e^{iw\sum_{n=N(t_{(m+m^{\prime})l+j^{\prime}-1})+1}^{N(t)}Z_{n}^{j^{\prime}}}\big).
\end{align*}
By similar method in the Proof of Lemma 2.2, we have
\begin{align*}
E\big(e^{iw\big(S_{t}-S_{s}\big)}\big)&=e^{iw\big(D_{t}-D_{s}\big)}e^{\int_{\mathbb{R}}
(e^{iwz}-1)\Big[\big(\Lambda(t_{j+ml})-\Lambda(s)\big)F_{j}
(dz)+\cdots+\big(\Lambda(t_{(m+1)l})-\Lambda(t_{(m+1)l-1})\big)F_{l}(dz)}\\
&\quad ^{+\sum_{k=1}^{m^{\prime}-1}\big(\Lambda(t_{(m+k)l+1})-\Lambda(t_{(m+k)l})\big)F_{1}
(dz)+\cdots+\sum_{k=1}^{m^{\prime}-1}\big(\Lambda(t_{(m+k+1)l})-\Lambda(t_{(m+k+1)l-1})\big)F_{l}(dz)}\\
&\quad ^{+\big(\Lambda(t_{(m+m^{\prime})l+1})-\Lambda(t_{(m+m^{\prime})l})F_{1}(dz)\big)+\cdots+
\big(\Lambda(t)-\Lambda(t_{(m+m^{\prime})l+j^{\prime}-1})F_{j^{\prime}}(dz)\big)\Big]}.
\end{align*}
Since $s+K\tau\in A_{j+(m+K)l}$ and $t+K\tau\in A_{j^{\prime}+(m+m^{\prime}+K)l},$ it follows the same method used in 
the computation of the characteristic function of $S_{t}-S_{s}$ that
\begin{align*}
E\big(e^{iw\big(S_{t+K\tau}-S_{s+K\tau}\big)}\big)=e^{iw\big(D_{t+K\tau}-D_{s+K\tau}\big)}& e^{\int_{\mathbb{R}}
(e^{iwz}-1)\Big[\big(\Lambda(t_{j+(m+K)l})-\Lambda(s+K\tau)\big)F_{j}(dz)+\cdots}\\
&^{+\big(\Lambda(t_{(m+K+1)l})-\Lambda(t_{(m+K+1)l-1})\big)F_{l}(dz)}\\
&^{+\sum_{k=1}^{m^{\prime}-1}\big(\Lambda(t_{(m+K+k)l+1})-\Lambda(t_{(m+K+k)l})\big)F_{1}(dz)+\cdots}\\
&^{+\sum_{k=1}^{m^{\prime}-1}\big(\Lambda(t_{(m+K+k+1)l})-\Lambda(t_{(m+K+k+1)l-1})\big)F_{l}(dz)}\\
&^{+\big(\Lambda(t_{(m+m^{\prime}+K)l+1})-\Lambda(t_{(m+m^{\prime}+K)l})F_{1}(dz)\big)+\cdots}\\
&^{+\big(\Lambda(t+K\tau)-\Lambda(t_{(m+m^{\prime}+K)l+j^{\prime}-1})F_{j^{\prime}}(dz)\big)\Big]}.
\end{align*}
By Definition 2.2 and Definition 2.3, we have for partition $0\leq t_{i}< t_{j}$ and $K\in\mathbb{N}$
\begin{align}\nonumber
\Lambda(t_{j})-\Lambda(t_{i})&=\Lambda(t_{j}+K\tau)-\Lambda(t_{i}+K\tau)\\ \label{equ6.1}
&=\Lambda(t_{j+Kl})-\Lambda(t_{i+Kl})
\end{align}
 and $D_{t}=D_{t+K\tau}$ for $t\geq0$. Thus
\begin{align*}
E\big(e^{iw\big(S_{t}-S_{s}\big)}\big)=E\big(e^{iw\big(S_{t+K\tau}-S_{s+K\tau}\big)}\big).
\end{align*}
\textbf{P3: Proof of Theorem 3.1}\\ 
Proof $a:$ Let $(S)_{t\geq0}$ be the simple semi-L\'evy process defined by (\ref{equ2.2}) and $Z_{i}$ be $i^{th}$ jump
size. Furthermore 
$T_{1}$ denote the time that the first jump occurs and $T_{j}, j=2, 3, \cdots,$ be the time intervals between the 
$(j-1)^{th}$ and $j^{th}$ jumps and $\Upsilon_{n}:=\sum_{j=1}^{n}T_{j},$ for $n\in\mathbb{N},$ and $\Upsilon_{0}=0$.
For $n\in\mathbb{N}$ 
\begin{align*}
Q_{n}&:=\big(I+(Z_{n})^{2}\mathbf{e}\mathbf{a}^{\prime}\big)e^{B(\Upsilon_{n}-\Upsilon_{n-1})},\\
R_{n}&:=\alpha_{0}(Z_{n})^{2}\mathbf{e}.
\end{align*}
It follows from \cite{b9} that $\mathbf{Y}_{t}$ satisfies in
\begin{align}\label{equ6.2}
\mathbf{Y}_{t}=\mathfrak{J}_{s,t}\mathbf{Y}_{s}+\mathfrak{K}_{s,t}, \qquad 0\leq s\leq t
\end{align}
where,
$$\mathfrak{J}_{s,t}=e^{B(t-\Upsilon_{N(t)})}Q_{N(t)}\cdots Q_{N(s)+2}\big(I+(Z_{N(s)+1})^{2}\mathbf{e}\mathbf{a}^{\prime}\big)e^{B(\Upsilon_{N(s)+1}-s)},$$
$$\mathfrak{K}_{s,t}=e^{B(t-\Upsilon_{N(t)})}\big(R_{N(t)}+Q_{N(t)}R_{N(t)-1}+\cdots+Q_{N(t)}\cdots Q_{N(s)+2}R_{N(s)+1}\big).$$
In order to prove that the sequence $\big(\mathfrak{J}_{s+k\tau,t+k\tau},\mathfrak{K}_{s+k\tau,t+k\tau}\big)_{k\geq0}$ is independently 
identically distributed, let $m\tau\leq s\leq t\leq (m+1)\tau$ such that $m\in\mathbb{N}^{0}.$ We define
\begin{align}\nonumber
\Upsilon_{1}^{(m+1)}&:=\Upsilon_{N(s)+1}-s,\quad Z_{1}^{(m+1)}:=Z_{N(s)+1},\\ \nonumber
&\vdots\\ \label{equ6.3}
\Upsilon_{N(t)-N(s)}^{(m+1)}&:=\Upsilon_{N(t)}-s, \quad Z_{N(t)-N(s)}^{(m+1)}:=Z_{N(t)}.
\end{align}
Therefore, $\big(\mathfrak{J}_{s,t},\mathfrak{K}_{s,t}\big)$ is function of random vector $\big(N(t)-N(s), \Upsilon^{(m+1)}_{1},\cdots, \Upsilon^{(m+1)}_{N(t)-N(s)}, Z^{(m+1)}_{1},$ $\cdots, Z^{(m+1)}_{N(t)-N(s)}\big).$  Using the fact that increments of the Poisson process $\lbrace N(t): t\geq0\rbrace$ are independent and the density function of random vector $(X_{1},\cdots,X_{n})$ can be computed as follows:
$$f_{\big(X_{1},\cdots,X_{n}\big)}(x_{1},\cdots,x_{n})=\lim_{\delta_{1},\cdots,\delta_{n}\rightarrow\infty}\frac{P\big(x_{1}<X_{1}\leq x_{1}+\delta_{1},\cdots,x_{n}<X_{n}\leq x_{n}+\delta_{n}\big)}{\delta_{1}\cdots\delta_{n}},$$
we can give the conditional density of $\big(\Upsilon_{1}^{(m+1)},\cdots, \Upsilon_{n}^{(m+1)}\vert N(t)-N(s)=n\big)$ as follows:
\begin{align*}
f_{\big(\Upsilon_{1}^{(m+1)},\cdots, \Upsilon_{n}^{(m+1)}\vert N(t)-N(s)=n\big)}(s_{1},\cdots,s_{n})&=\frac{n!
\lambda(s_{1})\times\cdots\times\lambda(s_{n})}{\big(\Lambda(t)-\Lambda(s)\big)^{n}},\\ 
&\qquad \qquad \qquad  0<s_{1}<\cdots< s_{n}<\tau.
\end{align*}
Since the increment process $N(t)-N(s)$ is a Poisson process with mean $\big(\Lambda(t)-\Lambda(s)\big)$ such that $\Lambda(t)-\Lambda(s)=\Lambda(t+k\tau)-\Lambda(s+k\tau),$ for all $k\in\mathbb{N}^{0},$ it follows from Definition 2.1 and conditional density that  
$$\big(\mathfrak{J}_{s,t},\mathfrak{K}_{s,t}\big)\overset{d}{=}\big(\mathfrak{J}_{s+k\tau,t+k\tau},\mathfrak{K}_{s+k\tau,t+k\tau}\big).$$
The independence of the sequence $\big(\mathfrak{J}_{s+k\tau,t+k\tau},\mathfrak{K}_{s+k\tau,t+k\tau}\big)$ is clear, since
$\mathfrak{J}_{s+k\tau,t+k\tau}$ and $\mathfrak{K}_{s+k\tau,t+k\tau}$ are constructed only from the segment $S_{u}, (s+k\tau)\leq u\leq
(t+k\tau),$ of the semi-Levy process $S$.\\
If $0\leq s\leq t,$ there are $n, m\in\mathbb{N}^{0}$ such that $s\in\big[n\tau,(n+1)\tau\big)$ and $t\in\big[(n+m)\tau,(n+m+1)\tau\big).$ By
iterating (\ref{equ3.1}) we obtain
\begin{align*}
\mathbf{Y}_{t}=&\Big(\mathfrak{J}_{(n+m)\tau,t}\mathfrak{J}_{(n+m-1)\tau,(n+m)\tau}\cdots\mathfrak{J}_{(n+1)\tau,
(n+2)\tau}\mathfrak{J}_{s,(n+1)\tau}\Big)\mathbf{Y}_{s}+\Big[\mathfrak{K}_{(n+m)\tau, t}\\
&+\mathfrak{J}_{(n+m)\tau, t}\mathfrak{K}_{(n+m-1)\tau, (n+m)\tau}+\cdots
+\mathfrak{J}_{(n+m)\tau, t}\mathfrak{J}_{(n+m-1)\tau, (n+m)\tau}\cdots
\mathfrak{J}_{(n+1)\tau,(n+2)\tau}\mathfrak{K}_{s,(n+1)\tau}\Big].
\end{align*}
It follows from \cite{b9} and (\ref{equ6.2}) that
$$\mathfrak{J}_{s,t}=\mathfrak{J}_{(n+m)\tau,t}\mathfrak{J}_{(n+m-1)\tau,(n+m)\tau}\cdots\mathfrak{J}_{(n+1)\tau,(n+2)\tau}\mathfrak{J}_{s,(n+1)\tau}$$
$$\mathfrak{K}_{s,t}=\mathfrak{K}_{(n+m)\tau, t}+\mathfrak{J}_{(n+m)\tau, t}\mathfrak{K}_{(n+m-1)\tau, (n+m)\tau}+\cdots+\mathfrak{J}_{(n+m)\tau, t}\cdots\mathfrak{J}_{(n+1)\tau,(n+2)\tau}\mathfrak{K}_{s,(n+1)\tau},$$
therefore, for all $k\in\mathbb{N}^{0},$
$$\big(\mathfrak{J}_{s,t},\mathfrak{K}_{s,t}\big)\overset{d}{=}\big(\mathfrak{J}_{s+k\tau,t+k\tau},\mathfrak{K}_{s+k\tau,t+k\tau}\big).$$

(b) By iterating (\ref{equ3.1}) we obtain
\begin{align*}
\mathbf{Y}_{t+m\tau}=\Big(\mathfrak{J}_{t+(m-1)\tau,t+m\tau}\cdots\mathfrak{J}_{t,t+\tau}\Big)\mathbf{Y}_{t}
+&\Big[\mathfrak{K}_{t+(m-1)\tau,t+m\tau}+\mathfrak{J}_{t+(m-1)\tau,t+m\tau}\mathfrak{K}_{t+(m-2)\tau,t+(m-1)\tau}\\
&+\cdots+\mathfrak{J}_{t+(m-1)\tau,t+m\tau}\cdots\mathfrak{J}_{t+\tau,t+2\tau}\mathfrak{K}_{t,t+\tau}\Big].
\end{align*}
Since $\big(\mathfrak{J}_{t+(m-1)\tau,t+m\tau},\mathfrak{K}_{t+(m-1)\tau,t+m\tau}\big)_{m\in\mathbb{N}}$, are independent and identically distributed it follows immediately
\begin{align*}
\mathbf{Y}_{t+m\tau}&\overset{d}{=}\prod_{k=1}^{m}\mathfrak{J}_{t+(k-1)\tau, t+k\tau}\mathbf{Y}_{t}
+\mathfrak{K}_{t,t+\tau}+\sum_{k=1}^{m-1}\mathfrak{J}_{t,t+\tau}\cdots\mathfrak{J}_{t+(k-1)\tau,t+k\tau}\mathfrak{K}_{t+k\tau, t+(k+1)\tau}.
\end{align*}
Note that the $\mathfrak{K}_{t,t+\tau}+\sum_{k=1}^{m-1}\mathfrak{J}_{t,t+\tau}\cdots\mathfrak{J}_{t+(k-1)\tau,t+k\tau}\mathfrak{K}_{t+k\tau, t+(k+1)\tau}$ is the partial sums of the infinite series
\begin{align}\label{equ6.4}
\mathbf{U}^{(t)}:=\mathfrak{K}_{t,t+\tau}+\sum_{k=1}^{\infty}\mathfrak{J}_{t,t+\tau}\cdots\mathfrak{J}_{t+(k-1)\tau,t+k\tau}\mathfrak{K}_{t+k\tau, t+(k+1)\tau}.
\end{align}
Thus, using the general theory of random recurrence equations (see Bougerol and Picard \cite{b6}, Brandt \cite{b7} and Vervaat \cite{v}) and condition (\ref{equ3.2}), we prove the almost sure absolute convergence of the series (\ref{equ6.4}). Let $P$ be such that $\Delta:=P^{-1}BP$ is diagonal. Then we have for $t\geq0,$
$$||e^{Bt}||_{B,r}=||Pe^{\Delta t}P^{-1}||_{B,r}=||P^{-1}Pe^{\Delta t}P^{-1}P||_{B,r}=||e^{\Delta t}||_{r}=e^{\eta t}.$$
Using (\ref{equ6.2}), (\ref{equ6.3}) and condition (\ref{equ3.2}) show
\begin{align*}
E\big(log||\mathfrak{J}_{t,t+\tau}||_{B,r}\big)&\leq E\Big(\eta\tau+log\big(1+(Z_{N(t+\tau)-N(\tau)}^{(2)})^{2}||
\mathbf{e}\mathbf{a}^{\prime}||_{B,r}\big)+\cdots+log\big(1+(Z_{1}^{(2)})^{2}||\mathbf{e}\mathbf{a}^{\prime}||_{B,r}\big)\\ 
&\quad+log\big(1+(Z_{N(\tau)-N(t)}^{(1)})^{2}||\mathbf{e}\mathbf{a}^{\prime}||_{B,r}\big)
+\cdots+log\big(1+(Z_{1}^{(1)})^{2}||\mathbf{e}\mathbf{a}^{\prime}||_{B,r}\big)\Big)<0\,
\end{align*}
and it follows from \cite{b9}, (\ref{equ6.2}) and (\ref{equ2.4})
\begin{align*}
E\big(log||\mathfrak{K}_{t,t+\tau}||_{B,r}\big)&\leq E\Big(\sum_{j=1}^{N(t+\tau)-N(t)}\Big[log(1+(Z_{N(t+\tau)-j+1})^{2}||ea^{\prime}||_{B,r})+(Z_{N(t+\tau)-j+1})^{2}\Big]\\
&\quad+log(\alpha_{0}||e||_{B,r})\Big)<\infty.
\end{align*}
Hence, the strong law of large numbers yield
$$\limsup_{k\rightarrow\infty}\dfrac{1}{k}\Big(\sum_{j=1}^{k}log||\mathfrak{J}_{t+(j-1)\tau,t+j\tau}||_{B,r}+log||\mathfrak{K}_{t+k\tau,t+(k+1)\tau}||_{B,r}\Big)<0\quad a.s.,$$
i.e.
$$\limsup_{k\rightarrow\infty}\Big(||\mathfrak{J}_{t,t+\tau}\cdots\mathfrak{J}_{t+(k-1)\tau,t+k\tau}\mathfrak{K}_{t+k\tau, t+(k+1)\tau}||_{B,r}\Big)^{\frac{1}{k}}<1\quad a.s.$$
From Cauchy's root criterion follows that series (\ref{equ6.4}) is almost sure absolute convergence. Since the state process $\mathbf{Y}$ has cadlag paths, it follows that $||\mathbf{Y}_{t}||_{B,r}$ is almost surely finite. Therefore
\begin{align*}
||\mathbf{Y}_{t+m\tau}-\mathbf{U}^{(t)}||_{B,r}\leq ||\prod_{k=1}^{m}\mathfrak{J}_{t+(k-1)\tau, t+k\tau}||_{B,r} \Big(||\mathbf{Y}_{t}-\mathbf{U}^{(t)}_{m}||_{B,r}\Big)\longrightarrow 0\quad a.s.
\end{align*}
where $\mathbf{U}^{(t)}_{m}:=\mathfrak{K}_{t+m\tau,t+(m+1)\tau}+\sum_{k=m}^{\infty}\mathfrak{J}_{t+m\tau,t+(m+1)\tau}
\cdots\mathfrak{J}_{t+k\tau,t+(k+1)\tau}\mathfrak{K}_{t+(k+1)\tau,t+(k+2)\tau}.$ It follows from \cite{b10} that $\mathbf{Y}_{t+m\tau}$
converges in distribution to $\mathbf{U}^{(t)},$ for fixed $t\in[0,\tau)$. That $\mathbf{U}^{(t)}$ satisfies (\ref{equ3.3}) and is the unique
solution is clear by the general theory of random recurrence equations.

($c$) It suffices to show that for any $s_{1}, s_{2}, \cdots,s_{n}$ and $k\in\mathbb{N}^{0}$
$$\big(\mathbf{Y}_{s_{1}}, \mathbf{Y}_{s_{2}}, \cdots, \mathbf{Y}_{s_{n}}\big)\overset{d}{=}\big(\mathbf{Y}_{s_{1}+k\tau}, \mathbf{Y}_{s_{2}+k\tau}, \cdots, \mathbf{Y}_{s_{n}+k\tau}\big).$$
Using the recursion equation (\ref{equ6.2}) and analysis is used in ($a$) we obtain above relation. We give the proof for $s_{1}\in[0,\tau)$ and $s_{2}\in[\tau,2\tau).$ The general case is similar. Therefore
$$\mathbf{Y}_{s_{2}}=\mathfrak{J}_{\tau,s_{2}}\mathfrak{J}_{s_{1},\tau}\mathbf{Y}_{s_{1}}+\mathfrak{K}_{\tau,s_{2}}
+\mathfrak{J}_{\tau,s_{2}}\mathfrak{K}_{s_{1},\tau},$$
and 
$$\mathbf{Y}_{s_{1}}=\mathfrak{J}_{0,s_{1}}\mathbf{Y}_{0}+\mathfrak{K}_{0,s_{1}}.$$
The random vector $\big(\mathbf{Y}_{s_{1}}, \mathbf{Y}_{s_{2}}\big)$ is function from $\big(\mathfrak{J}_{0,s_{1}}, \mathfrak{K}_{0,s_{1}}, \mathfrak{J}_{s_{1},\tau}, \mathfrak{K}_{s_{1},\tau}, \mathfrak{J}_{\tau,s_{2}}, \mathfrak{K}_{\tau,s_{2}}, \mathbf{Y}_{0}\big)$ and with similar argument also shows that the random vector $\big(\mathbf{Y}_{s_{1}+k\tau}, \mathbf{Y}_{s_{2}+k\tau}\big)$ is function from $\big(\mathfrak{J}_{k\tau,s_{1}+k\tau},$ $\mathfrak{K}_{k\tau,s_{1}+k\tau}, \mathfrak{J}_{s_{1}+k\tau,(k+1)\tau}, \mathfrak{K}_{s_{1}+k\tau,(k+1)\tau}, \mathfrak{J}_{(k+1)\tau,s_{2}+k\tau}, \mathfrak{K}_{(k+1)\tau,s_{2}+k\tau}, \mathbf{Y}_{k\tau}\big).$ Using ($a$) and assumption $\mathbf{Y}_{0}\overset{d}{=}\mathbf{U}^{(0)},$ it follows that $\big(\mathbf{Y}_{s_{1}}, \mathbf{Y}_{s_{2}}\big)\overset{d}{=}\big(\mathbf{Y}_{s_{1}+k\tau}, \mathbf{Y}_{s_{2}+k\tau}\big).$\\ \\
\textbf{P4: Proof of Corollary 3.2}\\
Since for all $s\in[t,t+p],$ the process $dS_{s}$ is independent of $\mathcal{F}_{s},$ it follows from (\ref{equ2.5}) and Corollary 3. that
\begin{align*}
E\big(G_{t}^{(p)}\big)&=\int_{t}^{t+p}E(\sqrt{V_{s}})E(S_{s+ds}-S_{s})\\
&\overset{d}{=}\int_{t}^{t+p}E(\sqrt{V_{s+k\tau}})E(S_{s+ds+k\tau}-S_{s+k\tau})=E\big(G_{t+k\tau}^{(p)}\big).
\end{align*}
In order to prove that the covariance function of $G_{t}^{(p)}$ is periodic, it suffices to show that 
$$E\big(G_{t}^{(p)}G_{t+h}^{(p)}\big)=E\big(G_{t+m\tau}^{(p)}G_{t+h+m\tau}^{(p)}\big).$$
Let $E_{t+p}$ denote conditional expectation with respect to the $\sigma$-algebra $\mathcal{F}_{t+p}.$ Since the increments of $S$ on the interval $(t+h,t+h+p]$ are independent of $\mathcal{F}_{t+p}$ and the increment process $G_{t}^{(p)}$ is measurable $\mathcal{F}_{t+p},$ we have
\begin{align*}
E\big(G_{t}^{(p)}G_{t+h}^{(p)}\big)&=E\Big(G_{t}^{(p)}E_{t+p}\big(G_{t+h}^{(p)}\big)\Big)\\
&=\int_{t+h}^{t+h+p}E\Big(\int_{t}^{t+p}\sqrt{V_{s}}\sqrt{V_{u}}dS_{u}\Big)E\big(dS_{s}\big).
\end{align*}
Since $\sqrt{V_{s}}\sqrt{V_{u}}dS_{u}$ is function of $\big(\mathfrak{J}_{u+du,s^{-}},\mathfrak{K}_{u+du,s^{-}},\mathfrak{J}_{u,u+du},
\mathfrak{K}_{u,u+du},\mathfrak{J}_{u^{-},u},\mathfrak{K}_{u^{-},u},\mathfrak{J}_{t,u^{-}},\mathfrak{K}_{t,u^{-}},$ 
$\mathbf{Y}_{t}\big)$
and this vector has the same distribution with $\big(\mathfrak{J}_{u+du+k\tau,s^{-}+k\tau},\mathfrak{K}_{u+du+k\tau,s^{-}+k\tau},$
$\mathfrak{J}_{u+k\tau,u+du+k\tau},$ $\mathfrak{K}_{u+k\tau,u+du+k\tau},\mathfrak{J}_{u^{-}+k\tau,u+k\tau},
\mathfrak{K}_{u^{-}+k\tau,u+k\tau},\mathfrak{J}_{t+k\tau,u^{-}+k\tau},\mathfrak{K}_{t+k\tau,u^{-}+k\tau},\mathbf{Y}_{t+k\tau}\big),$ it
 follows that 
$$cov\big(G_{t}^{(p)},G_{t+h}^{(p)}\big)=cov\big(G_{t+k\tau}^{(p)},G_{t+h+k\tau}^{(p)}\big).$$
\textbf{P5: Proof of Lemma 4.1}\\
($a$) Let $\tilde{\mathbf{Y}}_{t}$ be state process of semi Levy driven COGARCH(1,1) process. Then
\begin{align}\label{equ6.5}
\tilde{\mathbf{Y}}_{t}=\tilde{\mathfrak{J}}_{0,t}\tilde{\mathbf{Y}}_{0}+\tilde{\mathfrak{K}}_{0,t},
\end{align}
where
\begin{align*}
\tilde{\mathfrak{J}}_{0,t}&=exp\Big(\eta t+\sum_{i=1}^{N(t)}log(1+||\mathbf{e}\mathbf{a}^{\prime}||_{B,r}Z_{i}^{2})\Big),
\end{align*}
\begin{align*}
\tilde{\mathfrak{K}}_{0,t}&=\alpha_{0}||\mathbf{e}||_{B,r}\times exp\Big(\sum_{i=1}^{N(t)}log(1+||\mathbf{e}\mathbf{a}^{\prime}||_{B,r}Z_{i}^{2})\Big)\sum_{i=1}^{N(t)}Z_{i}^{2}.
\end{align*}
It follows from \cite{b9} that for all $t\geq0$
\begin{align}\label{equ6.6}
||\mathfrak{J}_{0,t}||_{B,r}\leq exp\Big(\eta t+\sum_{i=1}^{N(t)}log(1+||\mathbf{e}\mathbf{a}^{\prime}||_{B,r}Z_{i}^{2})\Big),
\end{align}
\begin{align}\label{equ6.7}
||\mathfrak{K}_{0,t}||_{B,r}\leq\alpha_{0}||\mathbf{e}||_{B,r}\times exp\Big(\sum_{i=1}^{N(t)}log(1+||\mathbf{e}\mathbf{a}^{\prime}||_{B,r}Z_{i}^{2})\Big)\sum_{i=1}^{N(t)}Z_{i}^{2}.
\end{align}
Now define a cadlag process $\lbrace X_{t}: t\geq0\rbrace$ by
$$X_{t}=-\eta t-\sum_{i=1}^{N(t)}log(1+||\mathbf{e}\mathbf{a}^{\prime}||_{B,r}Z_{i}^{2}),\quad t\geq0.$$
Then, $X_{t}$ is a negative simple pure jump semi Levy procress. It follows from Definition 2.1 and  Remark 3.2 that
\begin{align*}
E\big(e^{-cX_{t}}\big)&=E\Big(exp\Big(c\eta t+c\sum_{i=1}^{N(t)}log(1+||\mathbf{e}\mathbf{a}^{\prime}||_{B,r}Z_{i}^{2})\Big)\Big)\\
&=exp\Big(c\eta t+\int_{0}^{t}\int_{\mathbb{R}}\big((1+||\mathbf{e}\mathbf{a}^{\prime}||_{B,r}z^{2})^{c}-1\big)\nu_{s}(dz)ds\Big).
\end{align*}
Using a similar analysis is used in proof of Proposition 3.2 \cite{k1}, it follows from \cite{b9} that
\begin{align*}
\tilde{\mathfrak{K}}_{0,t}=||\mathbf{e}\mathbf{a}^{\prime}||_{B,r}^{-1}\alpha_{0}||\mathbf{e}||_{B,r}\Big[e^{-X_{t}}-\eta\int_{0}^{t}e^{-(X_{t}-X_{u})}du-1\Big].\end{align*}
It follows from (\ref{equ6.5}), (\ref{equ6.6}) and (\ref{equ6.7}) that $||\mathbf{Y}_{t}||_{B,r}<\tilde{\mathbf{Y}}_{t},$ for all $t\geq0$. Thus $E(S_{t}^{2})<\infty$ and $E(S_{t}^{4})<\infty$ imply $E(\mathbf{Y}_{t})<\infty$ and $cov(\mathbf{Y}_{t})<\infty,$ respectively.

In the proof of Theorem 3.1($a$) we have seen that (\ref{equ3.2}) implies that the sequence $\tilde{\mathbf{Y}}_{m\tau}$ converges in distribution to a finite random vector $\tilde{\mathbf{U}}$ which of the vector $\tilde{\mathbf{U}}$ is the unique solution of the random equation
$$\tilde{\mathbf{U}}\overset{d}{=}\tilde{\mathfrak{J}}_{0}\tilde{\mathbf{U}}+\tilde{\mathfrak{K}}_{0},$$
where $\big(\tilde{\mathfrak{J}}_{0},\tilde{\mathfrak{K}}_{0}\big)\overset{d}{=}\big(\tilde{\mathfrak{J}}_{0,\tau},\tilde{\mathfrak{K}}_{0,\tau}\big)$ and $\tilde{\mathbf{U}}$ is independent of $\big(\tilde{\mathfrak{J}}_{0},\tilde{\mathfrak{K}}_{0}\big).$
It follows from (\ref{equ6.4}), (\ref{equ6.6}) and (\ref{equ6.7}) that $\mathbf{U}\leq\tilde{\mathbf{U}}.$ Thus $E(S_{t}^{2})<\infty$ and $E(S_{t}^{4})<\infty$ imply $E(\mathbf{U})<\infty$ and $cov(\mathbf{U})<\infty,$ respectively.\\ \\
\textbf{P6: Proof of Lemma 4.2}\\
Using  (\ref{equ5.1}) and independence $\mathbf{Y}_{m\tau}$ and $(\mathfrak{J}_{m\tau, s_{1}+m\tau},\mathfrak{K}_{m\tau,s_{1}+m\tau}),$ we obtain 
\begin{align*}
E(\mathbf{Y}_{t})&=E(\mathfrak{J}_{m\tau,s_{1}+m\tau})E(\mathbf{Y}_{m\tau})+E(\mathfrak{K}_{m\tau,s_{1}+m\tau})\\
&=E(\mathfrak{J}_{0,s_{1}})E(\mathbf{U})+E(\mathfrak{K}_{0,s_{1}})
\end{align*}
where the last equality follows from that $(\mathfrak{J}_{m\tau,s_{1}+m\tau},\mathfrak{K}_{m\tau,s_{1}+m\tau})\overset{d}{=}(\mathfrak{J}_{0,s_{1}},\mathfrak{K}_{0,s_{1}})$ and assumption of section ($c$) of the Theorem 3.1. 

For computing  $cov(\mathbf{Y}_{t}, \mathbf{Y}_{t+h})$ it is sufficient to obtain $E(\mathbf{Y}_{t+h}\mathbf{Y}_{t}^{\prime}).$ It will therefore be followed from recursion equations which used in the proof of Theorem 4.1 that
\begin{align*}
\mathbf{Y}_{t}&=\mathfrak{J}_{m\tau,t}\mathbf{Y}_{m\tau}+\mathfrak{K}_{m\tau,t},\\
\mathbf{Y}_{t+h}&=\Big(\mathfrak{J}_{n\tau,t+h}\mathfrak{J}_{(n-1)\tau,n\tau}\cdots\mathfrak{J}_{m\tau,(m+1)\tau}\Big)\mathbf{Y}_{m\tau}\\
&\quad+\Big[\mathfrak{K}_{n\tau,t+h}+\mathfrak{J}_{n\tau,t+h}\mathfrak{K}_{(n-1)\tau,n\tau}+\sum_{i=1}^{n-m-1}\mathfrak{J}_{n\tau,t+h}\cdots\mathfrak{J}_{(n-i)\tau,(n-i+1)\tau}\mathfrak{K}_{(n-i-1)\tau,(n-i)\tau}\Big]
\end{align*}
The relation (\ref{equ4.2}) follows from independence the sequence $(\mathfrak{J}_{k\tau,s+k\tau},\mathfrak{K}_{k\tau,s+k\tau})$ for any $k\in\mathbb{N}^{0}$ and $s\in[0,\tau]$ and also independence $\mathbf{Y}_{m\tau}$ from this sequence for any $k\geq m.$\\ \\
\textbf{P7: Proof of Corollary 4.3}\\
Since for fixed $t$, almost surely $V_{t}=V_{t^{+}}=\alpha_{0}+\mathbf{a}^{\prime}\mathbf{Y}_{t},$ we have the expected value and covariance function volatility process from (\ref{equ2.6}).\\ \\
\textbf{P8: Proof of Proposition 4.4}\\
($a$) We imitate the proof of Theorem 6.1 of Brockwell, Chadraa, and Lindner \cite{b9}. Since $S$ is a martingale with zero mean, we have (\ref{equ4.3}). It follows from Ito isometry for square integrable martingales as integrators (e.g. \cite{r1}, IV 27) that
$$E\big(G_{t}^{(p)}G_{t+h}^{(p)}\big)=E\int_{0}^{t+h+p}V_{s}I_{[t,t+p)}(s)I_{[t+h,t+h+p}(s)d[S,S]_{s}=0,$$
and hence (\ref{equ4.4}) follows.\\
($b$) It follows from partial integration that
\begin{align}\nonumber
(G_{t}^{(p)})^{2}&=2\int_{t^{+}}^{t+p}G_{s^{-}}dG_{s}+[G,G]_{t^{+}}^{t+p}\\ \label{equ6.8}
&=2\int_{t}^{t+p}G_{s^{-}}\sqrt{V_{s}}dS_{s}+\sum_{t<s\leq t+p}V_{s}(\Delta S_{s})^{2}.
\end{align}
By similar analysis is used in ($a$), the compensation formula and \cite{c2} we have
\begin{align*}
E\big((G_{t}^{(p)})^{2}\big)=E\sum_{t<s\leq t+p}V_{s}(\Delta S_{s})^{2}=\int_{t}^{t+p}\int_{\mathbb{R}}E(V_{s})z^{2}\nu_{s}(dz)ds
\end{align*}
From Remark 3.2 the relation (\ref{equ4.5}) follows.\\
For proof of (\ref{equ4.6}), Since the increments of $S$ on the interval $(t,t+p]$ are independent of $\mathcal{F}_{t+p}$ and $S$ has expectation 0, it follows that
$$E_{t+p}\int_{t}^{t+p}G_{s^{-}}\sqrt{V_{s}}dS_{s}=0.$$
Thus it follows from the compensation formula and (\ref{equ6.2}) that
\begin{align*}
E_{t+p}\big((G_{t+h}^{(p)})^{2}\big)&=E_{t+p}\sum_{t+h<s\leq t+h+p}\big(\alpha_{0}+\mathbf{a}^{\prime}\mathfrak{J}_{t+p, s^{-}}\mathbf{Y}_{t+p}+a^{\prime}\mathfrak{K}_{t+p, s^{-}}\big)(\Delta S_{s})^{2}\\
&=\int_{t+h}^{t+h+p}\int_{\mathbb{R}}\big(\alpha_{0}+\mathbf{a}^{\prime}E(\mathfrak{J}_{t+p,s^{-}})\mathbf{Y}_{t+p}+
a^{\prime}E(\mathfrak{K}_{t+p,s^{-}})\big)z^{2}\nu_{s}(dz)ds,
\end{align*}
therefore
\begin{align*}
cov\big((G_{t}^{(p)})^{2},(G_{t+h}^{(p)})^{2}\big)=E\Big((G_{t}^{(p)})^{2}E_{t+p}\big((G_{t+h}^{(p)})^{2}\big)\Big)-E\big((G_{t}^{(p)})^{2}\big)E\big((G_{t+h}^{(p)})^{2}\big),
\end{align*}
and by Remark 3.2 and (\ref{equ2.6}) we have (\ref{equ4.6}). To calculate $cov\big(\mathbf{Y}_{t+p},(G_{t}^{(p)})^{2}\big),$ partial integration (\ref{equ6.8}) to get 
$$cov\big(\mathbf{Y}_{t+p},(G_{t}^{(p)})^{2}\big)=2cov\big(\mathbf{Y}_{t+p},\int_{t}^{t+p}G_{s^{-}}\sqrt{V_{s}}d\bar{S}_{s}\big)
+cov\big(\mathbf{Y}_{t+p},\int_{t^{+}}^{t+p}V_{s}d[S,S]_{s}\big).$$
To calculate the first term, let $I_{t}:=\int_{0}^{t}G_{s^{-}}\sqrt{V_{s}}dS_{s}.$ We know $E(I_{t})=0$ for all $t\geq0.$ Therefore
$$cov\big(\mathbf{Y}_{t+p},\int_{t}^{t+p}G_{s^{-}}\sqrt{V_{s}}dS_{s}\big)=E\big(I_{t+p}\mathbf{Y}_{t+p}\big)-E(\mathfrak{J}_{t,t+p})E\big(I_{t}\mathbf{Y}_{t}\big)-E(I_{t})E(\mathfrak{K}_{t,t+p}).$$
From \cite{b9}, partial integration and substituting $dV_{t^{+}}=a^{\prime}B\mathbf{Y}_{t}dt+\alpha_{0}V_{t}d[S,S]_{t}$ it follows that
\begin{align*}
E(I_{t}V_{t^{+}})&=E\int_{0}^{t}I_{s^{-}}dV_{s^{+}}+E\int_{0}^{t}V_{s}dI_{s}+E\big([V_{+},I]_{t}\big)\\
&=\mathbf{a}^{\prime}B\int_{0}^{t}E(I_{s^{-}}\mathbf{Y}_{s})ds+\alpha_{q}\int_{0}^{t}\int_{\mathbb{R}}E(I_{s^{-}}V_{s})z^{2}\nu_{s}(dz)ds\\
&\quad+E\int_{0}^{t}G_{s^{-}}\sqrt{V_{s}}V_{s}dS_{s}+\alpha_{q}E\int_{0}^{t}G_{s^{-}}\sqrt{V_{s}}V_{s}dM_{s},
\end{align*}
where $M_{s}:=\sum_{0<u\leq s}(\Delta S_{s})^{3}$ is a locally integrable martingale, with mean zero as a result of assumption that $\int_{\mathbb{R}}z^{3}\nu_{s}(dz)=0,$ for all $s\geq0.$ Thus, using the fact that $E\int_{0}^{t}G_{s^{-}}\sqrt{V_{s}}V_{s}dS_{s}=0,$ $E(I_{t}V_{t^{+}})=a^{\prime}E(I_{t}\mathbf{Y}_{t})$ and that $I_{s}\mathbf{Y}_{s}=I_{s^{-}}\mathbf{Y}_{s}=I_{s^{-}}\mathbf{Y}_{s^{-}}$ almost surely for fixed s, so we have
\begin{align*}
\mathbf{a}^{\prime}E(I_{t}\mathbf{Y}_{t})=\mathbf{a}^{\prime}B\int_{0}^{t}E(I_{s}\mathbf{Y}_{s})ds+\alpha_{0}\mathbf{a}^{\prime}
\int_{0}^{t}\int_{\mathbb{R}}E(I_{s}\mathbf{Y}_{s})z^{2}\nu_{s}(dz)ds,
\end{align*} 
The equality holds for any vector $\mathbf{a}$, hence
$$E(I_{t}\mathbf{Y}_{t})=B\int_{0}^{t}E(I_{s}\mathbf{Y}_{s})ds+\alpha_{0}\int_{0}^{t}\int_{\mathbb{R}}E(I_{s}\mathbf{Y}_{s})z^{2}\nu_{s}(dz)ds.$$\\
To calculate the second term of the covariance, it follows from \cite{b9} and (\ref{equ2.7}) that
\begin{align*}
cov\big(\mathbf{Y}_{t+p},\int_{t^{+}}^{t+p}V_{s}d[S,S]_{s}\big)&=cov\big(\mathbf{Y}_{t+p},(\mathbf{Y}_{t+p}^{\prime}-\mathbf{Y}_{t}^{\prime}-\int_{t^{+}}^{t+p}\mathbf{Y}_{s^{-}}dsB^{\prime})\mathbf{e}\big)\\
&=\big(cov(\mathbf{Y}_{t+p})-cov(\mathbf{Y}_{t+p},\mathbf{Y}_{t})-\int_{t^{+}}^{t+p}cov(\mathbf{Y_{t+p}},\mathbf{Y}_{s^{-}})dsB^{\prime}\big)\mathbf{e}.\end{align*}

\end{document}